\documentclass[12pt]{amsart}
\usepackage{amsmath,amssymb,amscd,array, mathrsfs }

\usepackage{amsmath,amscd,amssymb,amsthm,array}
\usepackage{color}

\usepackage{amsmath,amssymb,mathrsfs,amsthm, tikz-cd,mathrsfs}
\usepackage{comment}
\def\url#1{\expandafter\s

\tring\csname #1\endcsname}

\def\mmat #1,#2,#3,#4,{\text{\small\arraycolsep=3pt $
\begin{pmatrix}#1&#2\\#3&#4\end{pmatrix}$}}

\usepackage{hyperref}
\usepackage{capt-of}

\usepackage{multirow}
\usepackage[all]{xy}

\usepackage{comments}
\newComments\SBe{Said}{blue}
\newComments\SBo{Sofiane}{blue}
\newComments\AM{Nacer}{blue}
\newComments\DL{DL}{red}
\newComments\QEh{QEh}{blue}

\def\mmat #1,#2,#3,#4,{\text{\small\arraycolsep=3pt $
\begin{pmatrix}#1&#2\\#3&#4\end{pmatrix}$}}

\usepackage{lscape}
\usepackage{tikz-cd}
\usepackage{enumerate}

\usepackage{DLdef1}
\usepackage{multicol}

\def\mmat #1,#2,#3,#4,{\text{\small\arraycolsep=3pt $
\begin{pmatrix}#1&#2\\#3&#4\end{pmatrix}$}}

\renewcommand {\ssbegin}[2][*]
 {\refstepcounter{subsection}%
\if#1*
\addcontentsline{toc}{subsection}{\thesubsection.\hskip 1pc #2}%
\else
\addcontentsline{toc}{subsection}{\thesubsection.\hskip 1pc #2. #1}%
\fi
 \def \secno {\gdef \secno {}{\ssecfont
\thesubsection.\hskip 2ex}%
 }%
 \begin{#2}}

\renewcommand {\sssbegin}[2][*]
 {\refstepcounter{subsubsection}
\if#1*
\addcontentsline{toc}{subsubsection}{\thesubsubsection.\hskip 1pc #2}%
\else
\addcontentsline{toc}{subsubsection}{\thesubsubsection.\hskip 1pc #2. #1}
\fi
 \def \secno {\gdef \secno {}{\ssecfont \thesubsubsection.\hskip 2ex}%
 }%
 \begin{#2}}

\renewcommand {\parbegin}[2][*]
 {\refstepcounter{paragraph}
\if#1*
\addcontentsline{toc}{paragraph}{\theparagraph.\hskip 1pc #2}%
\else
\addcontentsline{toc}{paragraph}{\theparagraph.\hskip 1pc #2. #1}
\fi
 \def \secno {\gdef \secno {}{\ssecfont \theparagraph.\hskip 2ex}%
 }%
 \begin{#2}}

\newcommand {\ce}{{\text{CE}}}
\newcommand {\res}{{\rm{res}}}

\newcommand {\po}{{\mathfrak{po}}}

\newcommand{\del}{\partial}

\rmnameii{etr}{etr}
\rmnameii{evv}{ev}

\DeclareMathOperator{\K}{\mathbb{K}}

\newcommand {\w}{\omega}

\newcommand{\black}{\color{black}}
\newcommand{\dM}{\mathrm{d}}

\newcommand{\LR}{\rm{LR}}

\newcommand{\HH}{\mathrm{H}}
\newcommand{\PA}{\mathrm{PA}}
\newcommand{\obs}{\mathrm{obs}}

\begin{document}

\title[Cohomology of Restricted Poisson
algebras in characteristic 2]{Cohomology of Restricted Poisson
algebras in characteristic 2}

\author{Sofiane Bouarroudj}

\address{Division of Science and Mathematics, New York University Abu Dhabi, P.O. Box 129188, Abu Dhabi, United Arab Emirates.}
\email{sofiane.bouarroudj@nyu.edu}

\author{Quentin Ehret}
\address {Division of Science and Mathematics, New York University Abu Dhabi, P.O. Box 129188, Abu Dhabi, United Arab Emirates.}
\email{qe209@nyu.edu}

\author{Jiefeng Liu}

\address{School of Mathematics and Statistics, Northeast Normal University, Changchun 130024, China}
\email{liujf534@nenu.edu.cn}

\thanks{S.B. and Q.E. were supported by the grant NYUAD-065, and J. Liu was supported by  NSFC (12371029, W2412041).}



\begin{abstract}
 In this paper, we  study restricted Poisson algebras in characteristic 2 and their relationship with restricted Lie-Rinehart algebras, for which we develop a cohomology theory and investigate abelian extensions. We also construct a full cohomology complex for restricted Poisson algebras in characteristic 2 that captures formal deformations and prove that it is isomorphic to the cohomology complex of a suitable restricted Lie-Rinehart algebra, under certain assumptions.
 A number of examples are provided in order to illustrate our constructions.
\end{abstract}


\maketitle

\thispagestyle{empty}
\setcounter{tocdepth}{2}
\tableofcontents

\section{Introduction}

\subsection{Poisson and Lie-Rinehart algebras} Poisson algebras emerged in the 1970s as an algebraic interpretation of the concept of Poisson brackets, which are well-known to physicists working in classical mechanics. They have a wide range of applications and appear in various contexts, such as integrable systems, Hamiltonian mechanics, string theory, quantum mechanics, differential geometry, representation theory, and quantization, to name a few. The Poisson structure plays a crucial role in solving problems related to the aforementioned topics. A very good reference covering both the geometric and algebraic aspects of Poisson algebras is the book by Laurent-Gengoux, Pichereau and Vanhaecke (\cite{LPV}). Roughly speaking, a Poisson algebra is a vector space endowed with an associative commutative multiplication and a Lie bracket that is a derivation of the associative product, see Eq. \eqref{eq:poisson-definition}. 

Poisson cohomology of a Poisson manifold was introduced by Lichnerowicz in \cite{Li}, while homology for a Poisson manifold was introduced by Koszul in \cite{Kos}. From an algebraic point of view, the Poisson cohomology of a Poisson algebra over a field of characteristic zero is deduced from the usual Chevalley-Eilenberg cohomology of Lie algebras by replacing the usual cochains space of skewsymmetric multilinear maps by skewsymmetric multilinear \textit{derivations} of the associative product of the Poisson algebra (see e.g. \cite[Section 4.1.2]{LPV}). Moreover, double Poisson algebras were introduced in \cite{VDB} in order to develop Poisson geometry in a noncommutative setting.

Poisson algebras are closely related to Lie-Rinehart algebras, which are an algebraic generalization of both Lie algebras and Lie algebroids. A Lie-Rinehart algebra is a triple $(A,L,\theta)$, where $A$ is an  associative commutative algebra, $L$ is a Lie algebra and an $A$-module, and $\theta:L\rightarrow\Der(A)$ is a morphism of Lie algebras and of $A$-modules, called \textit{anchor}, which satisfies a compatibility condition that roughly controls the lack of $A$-linearity of the Lie bracket of $L$. They appeared in the works of Herz (\cite{He}), Palais (\cite{P}) and most notably Rinehart  (\cite{Ri}), who used this notion to develop a formalism of differential forms for general commutative algebras and also developed the representations and cohomology of Lie-Rinehart algebras.  More recently, Huebschmann explored the connection between Lie-Rinehart algebras and Poisson algebras in \cite{Hu}. He constructed a Lie-Rinehart algebra from a Poisson algebra using its module of Kähler differentials, in order to define Poisson cohomology and homology in terms of the corresponding structures of this Lie-Rinehart algebra. He also applied this construction to an arbitrary Poisson manifold and introduced the notion of smooth Poisson cohomology which is related to Lichnerowicz's work, and constructed an isomorphism between the algebraic and the smooth cohomology.

\subsection{Restricted Lie, Poisson and Lie-Rinehart algebras}
Restricted Lie algebras were introduced over a field of positive characteristic $p\geq 2$ by Jacobson in \cite{J} to give a `good' algebraic interpretation of properties of the derivations algebra of an associative algebra. They appear naturally in integration of algebraic groups, in representation theory and are useful in classification problems (see e.g \cite{SF}). Those algebras are equipped with a so-called $p$-map which resembles the $p$-th power of an associative algebra, see Section \ref{sec:resLie} for the full definition. The cohomology of restricted algebra is much more complicated than the usual Chevalley-Eilenberg cohomology in characteristic zero. It was first defined by Hochschild (see \cite{Hoch}), while explicit formulas for the low-degree cohomology appeared much later in the work of Evans and Fuchs (see \cite{EF}), who also studied central extensions of restricted Lie algebras. Interestingly, some of these formulas had already appeared earlier in a lesser-known paper by Pareigis \cite{Pa}, dating from 1967.\footnote{This paper, written in German, seems to have been largely forgotten, but it is of considerable interest to anyone studying restricted cohomology.} Another approach to the restricted cohomology consists in using the notion of Beck modules and a Quillen-Barr-Beck construction to define restricted cohomology groups and has been done by Dokas (see \cite{Do1}). More recently, formal deformations of restricted Lie algebras have been investigated in \cite{EM}.

Restricted Lie-Rinehart were formally introduced by Dokas \cite{Do2}, but they already appeared as early as 1955 in Hochschild's works (see \cite{Ho}) and in \cite{Ru}. In comparison to (ordinary) Lie-Rinehart algebras, the underlying Lie algebra is assumed to be restricted, as well as the anchor map, and there is an extra compatibility condition that roughly controls the lack of $p$-homogeneity of the $p$-map with respect to the action of the underlying associative algebra, see \cite[Definition 1.7]{Do2}. A cohomology theory for those objects based on Quillen-Barr-Beck methods has been introduced by Dokas in \cite{Do2}, while another approach adapted to formal deformations is proposed in \cite{Eh}.

The notion of restricted Poisson algebra was introduced by Bezrukavnikov and Kaledin (see \cite{BK2}) in order to investigate deformation quantization of algebraic manifolds in positive characteristics, as a continuation to their previous work on Fedosov quantization in characteristic zero (see \cite{BK1}). An equivalent definition was given by Bao, Ye and Zhang in \cite{BYZ}. They studied the relationship between Poisson algebras and restricted Lie-Rinehart algebras over field of characteristic $p\geq 3$. In particular, they showed that any restricted Poisson algebra gives rise to a restricted Lie-Rinehart algebra on the module of K\"{a}hler differentials, see \cite[Theorem 8.2]{BYZ}. Over a field of characteristic $p=2$, the definition of a restricted Poisson algebra was given by Petrogradsky and Shestakov in \cite{PS}. 

\subsection{The characteristic $p=2$ case. Outline of the paper and main results}

All of the studies mentioned above were conducted for $p>2$. The aim of this paper is to complete the picture for $p=2$. In that case, the restricted cohomology is known at any order (see \cite{EM}). Moreover, the equivalence between Bezrukavnikov-Kaledin's (see \cite{BK2}) and Bao-Ye-Zhang's (see \cite{BYZ}) definitions of restricted Poisson algebras is no longer valid and a new approach is needed.

In Section \ref{sec:background}, we review basics on restricted Lie algebras, their cohomology in characteristic $2$ and restricted Lie-Rinehart algebras. In Section \ref{sec:ab-ext-RLR}, we investigate abelian extensions of restricted Lie-Rinehart algebras in characteristic $2$. We introduce a well-defined cohomology for those objects derived from the restricted cohomology (see Proposition \ref{prop:RLRcoho}) and show that under some conditions, the abelian extensions of restricted Lie-Rinehart algebras are classified up to equivalence by the second cohomology group (see Theorem \ref{thmext}). We conclude this section with a brief discussion comparing our work with that of Dokas (\cite{Do2}). Section \ref{sec:poisson} is devoted to the study of restricted Poisson algebras in characteristic 
$p=2$. We show that Petrogradsky and Shestakov's definition (see \cite{PS}) 
extends those in \cite{BK2} and \cite{BYZ} to characteristic 2. We also present some classical constructions, including morphisms, modules, and derivations, and provide a class of examples arising from the deformation theory of associative algebras (see Section \ref{sec:exdefoquant}). Furthermore, we introduce a cohomology for restricted Poisson algebras for $p=2$, see Proposition \ref{prop:cohopoisson}, and explore its relationship with the cohomology of restricted Lie–Rinehart algebras. Analogous to Huebschmann's construction in characteristic zero (see \cite{Hu}) and that of Bao–Ye–Zhang in characteristic $p \geq 3$ (see \cite{BYZ}), we show that any restricted Poisson algebra gives rise to a restricted Lie–Rinehart structure on its module of K\"{a}hler differentials (Theorem \ref{thm:RLRfromPoisson}). Moreover, when this module is free, we prove that the two cohomologies are isomorphic (see Theorem \ref{thm:coho-iso}). In Section \ref{sec:defopoisson}, we consider formal deformations of restricted Poisson algebras, and show that equivalent infinitesimal deformations are controlled by the second cohomology space (see Corollary \ref{cor:equiv-defo}) and study the obstructions to extending a given deformation to higher orders, see Proposition \ref{prop:obs}. Finally, we compute some examples of restricted Poisson algebras and their cohomology in Section \ref{sec:ex}.\\

\noindent\textbf{Conventions and notations.} Throughout the paper, `ordinary" shall be understood as ``not restricted" and unless otherwise specified (briefly in Section \ref{sec:resLie}), $\K$ denotes a field of characteristic $p=2$.






\section{Background}\label{sec:background}

\subsection{Restricted Lie algebras}\label{sec:resLie}
Let $L$ be a finite-dimensional Lie algebra over a field $\K$ of positive characteristic $p>0$ . Following \cite{J}, a map $(\cdot)^{[p]}:L\rightarrow L, \quad x\mapsto x^{[p]}$ is called a~\textit{$p$-map} (or $p$\textit{-structure}) on $L$ and $L$ is said to be {\it restricted}  if 
\begin{align}\label{def:restrictedLie}
(\lambda x)^{[p]}&=\lambda^p x^{[p]} \text{ for all } x\in L \text{ and for all } \lambda \in \K;\\
\ad_{x^{[p]}}&=(\ad_x)^p \text{ for all } x\in L;\\
(x+y)^{[p]}&=x^{[p]}+y^{[p]}+\displaystyle\sum_{1\leq i\leq p-1}s_i(x,y), \text{ for all } x,y\in L;\label{eq:si}
\end{align}
where the coefficients $s_i(x,y)$ can be obtained from the expansion
\begin{equation}
(\ad_{\lambda x+y})^{p-1}(x)=\displaystyle\sum_{1\leq i \leq p-1} is_i(x,y) \lambda^{i-1}.
\end{equation}
In the case where $p=2$, Condition \eqref{eq:si} reduces to
$$(x+y)^{[2]}=x^{[2]}+y^{[2]}+[x,y],~\forall x,y\in L.$$

\sssbegin{Example}  Any associative algebra $A$ can be turned into a restricted Lie algebra with the bracket $[a,b]=ab-ba,~\forall a,b\in A$ and the $2$-map $a^{[2]}=a^2,~\forall a\in A$, see \cite{J}. 
\end{Example}

Let $\left( L,[\cdot,\cdot]_L,(\cdot)^{[p]_L}\right) $ and $\left( H,[\cdot,\cdot]_H,(\cdot)^{[p]_H}\right) $ be two restricted Lie algebras. A Lie algebra morphism $\varphi:L\rightarrow H$ is called \textit{restricted morphism} (or \textit{$p$-morphism}) if it satisfies \begin{equation}
\varphi\bigl(x^{[p]_L} \bigl)=\varphi(x)^{[p]_H}\quad\forall x\in L.
\end{equation}
An $L$-module $M$ is called \textit{restricted} if in addition to the ordinary condition
 \begin{equation}
   [x,y]_L\cdot m=x\cdot(y\cdot m)-y\cdot(x\cdot m) ,~\forall x,y\in L,\forall m\in M,
\end{equation} it also satisfies
\begin{equation}
\underbrace{x\cdots x}_{p\text{~~terms}}\cdot m =x^{[p]_L}\cdot m, \quad  \forall x\in L,~\forall  m\in M.
\end{equation}
A linear map $d:L\rightarrow L$ is called a \textit{restricted derivation} of $L$  if in addition to the ordinary condition
 \begin{equation}
    d\bigl([x,y]_L\bigl)=[d(x),y]_L+[x,d(y)]_L,~\forall x,y\in L,
\end{equation} it also satisfies
\begin{equation}
    d\bigl(x^{[p]_L}\bigl)=\ad_x^{p-1}\circ d(x),~\forall x\in L.
\end{equation}

The following Theorem, due to Jacobson (see \cite{J}) is useful to study $p$-structures on Lie algebras.
\sssbegin{Theorem}\label{thm:jacobson}
Let $(L,[-,-])$ be a Lie algebra and let $(e_j)_{j\in J}$ be a~basis of $L$ such that there are $f_j\in L$ satisfying $(\ad_{e_j})^p=\ad_{f_j}$. Then, there exists exactly one $p$-map  $(\cdot)^{[p]}:L\rightarrow L$ such that 
\[
e_j^{[p]}=f_j \quad \text{ for all $j\in J$}.
\]
\end{Theorem}

\subsection{Restricted cohomology for $p=2$}\label{sec:rescoho2} In this Section, we recall the construction of the restricted cohomology for restricted Lie algebras over a field $\K$ of characteristic $p=2$. The complex was described in \cite{EM}, although it had been known among experts for some time.  

Let us first start with the Chevalley-Eilenberg cohomology.

\subsubsection{Chevalley-Eilenberg cohomology}	
Let $L$ be a Lie algebra and $M$ be a $L$-module.  The Chevalley-Eilenberg cochains of $L$ with coefficients in $M$ are defined by
\begin{align*}
C^m_{\rm CE}(L,M)&=\Hom_{\K}(\Lambda^mL, M)~\text{ for }m \geq 1,\\
C^0_{\rm CE}(L,M)&\cong M.
\end{align*}
The coboundary operators  ${\rm d}^m_{\rm CE}: C^m_{\rm CE}(L,M)\longrightarrow C^{m+1}_{\rm CE}(L,M)$ are defined by
\begin{align}\label{eq:CE-coboundary operator}
\nonumber {\rm d}_{\rm CE}^m(\varphi)(x_1,\cdots,x_{m+1})&=\sum_{1\leq i<j\leq m+1}(-1)^{i+j-1}\varphi\left([x_i,x_j],x_1,\cdots\hat{x_i},\cdots,\hat{x_j},\cdots, x_{m+1} \right)\\&+\sum_{i=1}^{m+1}(-1)^{i}x_i\varphi(x_1,\cdots,\hat{x_i},\cdots,x_{m+1}). 
 \end{align}
 
We denote  the cohomology groups of Chevalley-Eilenberg of $L$ with coefficients in $M$ by
$  \mathrm{H}^m_{\rm CE}(L,M).   $

\subsubsection{Restricted cochains for $p=2$}
Let $\bigl(L,[\cdot,\cdot],(\cdot)^{[2]}\bigl)$ be a restricted Lie algebra and let $M$ be a restricted $L$-module. We set $C_{\mathrm{res}}^0(L;M):=C_{\mathrm{CE}}^0 (L;M)$ and $C_{\res}^1(L;M):=C_{\mathrm{CE}}^1 (L;M)$.

Let $n\geq2$, $\varphi\in C^n_{\ce} (L,M)$, $\omega:L\otimes \wedge^{n-2} L\rightarrow M$, $\lambda\in\K$ and $x,z_2,\cdots,z_{n-1}\in L$. The pair $(\varphi,\omega)$ is a $n$-cochain of the restricted cohomology if
		\begin{align}
			\omega(\lambda x, z_2,\cdots,z_{n-1})=&~\lambda^2\omega(x,z_2,\cdots,z_{n-1})\\
			\omega(x+y,z_2,\cdots,z_{n-1})=&~\omega(x,z_2,\cdots,z_{n-1})+\omega(y,z_2,\cdots,z_{n-1})\\\nonumber&~+\varphi(x,y,z_2,\cdots,z_{n-1}),\\
        (z_2,\cdots, z_{n-1})\mapsto &~\w(\cdot,z_2,\cdots, z_{n-1}) \text{ is linear}.
        \end{align}
        
		We denote by $C_{\res}^n(L,M)$ the space of $n$-cochains of $L$ with values in $M$.

\subsubsection{Restricted coboundary operators for  $p=2$}\label{sec:diff.op.p=2}
For $n\geq 2$, the coboundary maps
	$\dM^n_{\res}:C^n_{\res}(L;M)\rightarrow C^{n+1}_{\res}(L;M)$ are given by $\dM^n_{\res}(\varphi,\omega)=\bigl(\dM^n_{\mathrm{CE}}(\varphi),\delta^n(\omega)\bigl),$ where
	\begin{align*}	\delta^n\omega(x,z_2,\cdots,z_n)&:=x\cdot\varphi(x,z_2,\cdots,z_n)+\sum_{i=2}^{n}z_i\cdot\omega(x,z_2,\cdots,\hat{z_i},\cdots,z_n)\\	&+\varphi(x^{[2]},z_2,\cdots,z_n)+\sum_{i=2}^{n}\varphi\left([x,z_i],x,z_2,\cdots,\hat{z_i},\cdots,z_n \right)\\
		&+\sum_{2\leq i<j\leq n}\omega\left(x,[z_i,z_j],z_2,\cdots,\hat{z_i},\cdots,\hat{z_j},\cdots,z_n  \right).
	\end{align*}
    For $n=0,1$, we define $\dM^0_{\res}=\dM^0_{\mathrm{CE}}$ and 
	\begin{align*}
		\dM^1_{\res}:C_{\res}^1(L;M)&\rightarrow C_{\res}^2(L,M)\\
		\varphi&\mapsto\bigl(\dM_{\mathrm{CE}}^1\varphi,\delta^1\varphi\bigl),~\text{where }\delta^1\varphi(x):=\varphi\bigl( x^{[2]} \bigl)+x\cdot\varphi(x),~\forall x\in L.
	\end{align*}

The complex $\left(C_{\res}^n(L;M),d^n_{\res} \right)_{n\geq 0}$ is a cochain complex, see, e.g., \cite[Theorem 4.9]{EM}. The $n^{th}$ restricted cohomology group of the Lie algebra $L$ in characteristic 2 is defined by	
		$$ \HH_{\res}^n(L;M):=Z_{\res}^n(L;M)/B^n_{\res}(L;M),$$ 
		with $Z_{\res}^n(L;M)=\Ker(\dM^n_{\res})$ the restricted $n$-cocycles and $B_{\res}^n(L;M)=\text{Im}(\dM^{n-1}_{\res})$ the restricted $n$-coboundaries. 

    In the case where $p\geq 3$, the restricted cohomology has been described in \cite{EF}.
        
\subsection{Restricted Lie-Rinehart algebras for $p=2$} (see \cite{Do2})\label{sec:RLRdef} A restricted Lie-Rinehart algebra is a triple $(A,L,\theta)$, where $A$ is an associative commutative algebra, $\bigl(L,[\cdot,\cdot],(\cdot)^{[2]}\bigl)$ is a restricted Lie algebra that is also an $A$-module and $\theta:L\rightarrow\Der(A)$ is an $A$-linear restricted Lie algebras morphism satisfying for all $x,y\in L$ and for all $a\in A$
\begin{align}\label{eq:leibniz}
    [x,a\cdot y]&=a\cdot [x,y]+\theta(x)(a)\cdot y;~ \text{ and }\\\label{eq:hochschild}
    (a\cdot x)^{[2]}&=a^2\cdot x^{[2]}+\theta(ax)(a)\cdot x.
\end{align}

Let $(A,L,\theta_L)$ and $(A,H,\theta_H)$ be restricted Lie-Rinehart algebras sharing the same associative commutative algebra $A$. A restricted Lie algebras morphism $\varphi:L\rightarrow H$ is called \textit{restricted Lie-Rinehart algebras morphism} if it is $A$-linear and if it satisfies $\theta_H\circ\varphi=\theta_L$.

A representation of a restricted Lie-Rinehart algebra $(A,L,\theta_L)$ is an $A$-module $M$ together with a restricted morphism of restricted Lie algebras $\rho:L\rightarrow \End(M)$ satisfying
\begin{equation}\label{eq:RLRmodule}
    \rho(x)(am)=a\rho(x)(m)+\theta_L(x)(a)m,~\forall x\in L,~\forall a\in A,~\forall m\in M.
\end{equation} Such a pair $(M,\rho)$ is also called \textit{restricted Lie-Rinehart module}.

\sssbegin{Example}
    Let $A$ be an associative commutative algebra and let $\Der(A)$ be its algebra of derivations. Then, $\Der(A)$ is a restricted Lie algebra with the bracket $[D_1,D_2]:=D_1\circ D_2-D_2\circ D_1,~ \forall D_1,D_2\in \Der(A)$ and the 2-map gievn by $D^{[2]}:=D\circ D,~ \forall D\in \Der(A)$. Moreover, $\Der(A)$ is an $A$-module with $(a\cdot D)(b):=aD(b),~\forall a,b\in A$ and $\forall D\in\Der(A)$. As a result, we obtain a restricted Lie-Rinehart algebra $(A,\Der(A),\id)$.
\end{Example}

The following result, due to Dokas, is a generalization of Jacobson's theorem \ref{thm:jacobson} to Lie-Rinehart algebras.
\sssbegin{Proposition}\cite[Proposition 2.2]{Do2}.\label{prop:jacobson-rinehart} Let $(A,L,\theta)$ be a Lie-Rinehart algebra such that $L$ is free as an $A$-module. Let $\{u_i\}_i$ be an ordered $A$-basis of $L$ such that there are elements $v_j\in L$ satisfying $\ad_{u_i}^p=\ad_{v_j}$ for all $i$. Then, the map $u_i\mapsto v_j$ defines a $p$-map on $L$ and $(A,L,\theta)$ is a restricted Lie-Rinehart algebra.
\end{Proposition}

\section{Abelian extensions of restricted Lie-Rinehart algebras, $p=2$}\label{sec:ab-ext-RLR}
\subsection{Cohomology} Let $(A,L,\theta)$ be a restricted Lie-Rinehart algebra in characteristic 2 and let $(M,\rho)$ be a restricted Lie-Rinehart module. For all $k\geq 1$, we denote by $\Hom_A(\wedge^kL,M)$ the space of $A$-linear maps $\varphi:\wedge^kL\rightarrow M$. For $n=0,1,$ we set $C_{\rm LR}^0(L;M):=M$ and $C_{\rm LR}^1(L;M):=\Hom_A(\wedge^1L,M)$. 

For $n\geq 2$ a pair $(\varphi,\w)\in C^n_{\rm res}(L;M)\cap\Hom_A(\wedge^nL,M)$ is a restricted Lie-Rinehart $n$-cochain if, for all $x,z_2,\cdots,z_{n-1}\in L$ and for all $a\in A$, we have
\begin{align}
\label{eq:res A-linear1}		\omega(a x, z_2,\cdots,z_{n-1})&=a^2\omega(x,z_2,\cdots,z_{n-1}),\\
\label{eq:res A-linear2}			\omega( x, z_2,\cdots,az_i,\cdots,z_{n-1})&=a\omega(x,z_2,\cdots,z_i,\cdots,z_{n-1}),
\end{align}

The space of all pairs $(\varphi,\w)\in C_{\res}^n(L;M)\cap\Hom_A(\wedge^nL,M)$ satisfying \eqref{eq:res A-linear1} and \eqref{eq:res A-linear2} is denoted by $C^n_{\rm LR}(L;M)$. The coboundary operators
${\rm d}_{\rm LR}^n:C^n_{\rm LR}(L;M)\rightarrow C^{n+1}_{\rm LR}(L;M)$ for
$n\geq 0$ are induced by the restricted coboundary operators for restricted Lie algebras, see Section \ref{sec:diff.op.p=2}.

\sssbegin{Proposition}\label{prop:RLRcoho}
	Let $(A,L,\theta)$ be a restricted Lie-Rinehart algebra, let $(M,\rho)$ be a restricted Lie-Rinehart module and let $(\varphi,\omega)\in C^{n}_{\LR}(L;M)$. Then, $\dM^{n}_{\rm LR}(\varphi,\omega)\in  C^{n+1}_{\LR}(L;M)$. Furthermore, we have $\dM^{n+1}_{\rm LR}\circ \dM^n_{\rm LR}=0$.
\end{Proposition}
Thus, we obtain a well-defined cochain complex $\bigl(\bigoplus_{n\geq0}C^{n}_{\rm LR}(L,M),\dM^n_{\rm LR}\bigl)$ and denote the $n$-th cohomology group of $\bigl(\bigoplus_{n\geq0}C^{n}_{\rm LR}(L,M),\dM^n_{\rm LR}\bigl)$ by $\HH_{\rm LR}^n(L,M)$.
\begin{proof}

Let $a\in A$ and $x,z_2,\cdots,z_{n}\in L$.	
We check that $\delta^n(\omega)$ satisfies \eqref{eq:res A-linear1} and \eqref{eq:res A-linear2}. Since $\omega$ satisfies \eqref{eq:res A-linear1}  and $\varphi$ is $A$-linear, we have 
\begin{eqnarray*}
& &\delta^n\omega(ax,z_2,\cdots,z_n)\\	
&=&\rho(ax)\varphi(ax,z_2,\cdots,z_n)+\sum_{i=2}^{n}\rho(z_i)\omega(ax,z_2,\cdots,\hat{z_i},\cdots,z_n)+\varphi((ax)^{[2]},z_2,\cdots,z_n)\\	&&+\sum_{i=2}^{n}\varphi\left([ax,z_i]_E,ax,z_2,\cdots,\hat{z_i},\cdots,z_n \right)+\sum_{1\leq i<j\leq n}\omega\left(ax,[z_i,z_j]_E,z_2,\cdots,\hat{z_i},\cdots,\hat{z_j},\cdots,z_n  \right)\\
	& &\\
	&=&a^2\rho(x)\varphi(x,z_2,\cdots,z_n)+a\theta_E(x)(a)\varphi(x,z_2,\cdots,z_n)\\
	&&+\sum_{i=2}^{n}a^2\rho(z_i)\omega(x,z_2,\cdots,\hat{z_i},\cdots,z_n)+\sum_{i=2}^{n}\theta_E(z_i)(a^2)\omega(x,z_2,\cdots,\hat{z_i},\cdots,z_n)\\
	&&+a^2\varphi(x^{[2]},z_2,\cdots,z_n)+a\theta_E(x)(a)\varphi(x,z_2,\cdots,z_n)\\
	&&+\sum_{i=2}^{n}a^2\varphi\left([x,z_i]_E,x,z_2,\cdots,\hat{z_i},\cdots,z_n \right)-\sum_{i=2}^{n}a\theta_E(z_i)(a)\varphi\left(x,x,z_2,\cdots,\hat{z_i},\cdots,z_n \right)\\
	&&+\sum_{1\leq i<j\leq n}a^2\omega\left(x,[z_i,z_j]_E,z_2,\cdots,\hat{z_i},\cdots,\hat{z_j},\cdots,z_n  \right)\\
	&=&a^2	\delta^n\omega(x,z_2,\cdots,z_n)+2a\theta_E(x)(a)\varphi(x,z_2,\cdots,z_n)\\
	&&+\sum_{i=2}^{n}2a\theta_E(z_i)(a)\omega(x,z_2,\cdots,\hat{z_i},\cdots,z_n)\\
	&=&a^2	\delta^n\omega(x,z_2,\cdots,z_n).
\end{eqnarray*}	

Furthermore, since $\omega$ satisfies \eqref{eq:res A-linear2}  and $\varphi$ is $A$-linear, we have 
\begin{eqnarray*}   
&&\delta^n\omega(x,az_2,\cdots,z_n)\\
&=&\rho(x)\varphi(x,az_2,\cdots,z_n)+\sum_{i=3}^{n}\rho(z_i)\omega(x,az_2,\cdots,\hat{z_i},\cdots,z_n)\\	
&&+\rho(az_2)\omega(x,z_3,\cdots,\hat{z_i},\cdots,z_n)+\sum_{j=3}^n\omega\left(x,[az_2,z_j]_E,z_3,\cdots,\hat{z_j},\cdots,z_n  \right)\\
&&+\sum_{i=3}^{n}\varphi\left([x,z_i]_E,x,az_2,\cdots,\hat{z_i},\cdots,z_n \right)+\varphi\left([x,az_2]_E,x,z_3,\cdots,z_n \right)\\
&&+\sum_{3\leq i<j\leq n}\omega\left(x,[z_i,z_j]_E,az_2,\cdots,\hat{z_i},\cdots,\hat{z_j},\cdots,z_n  \right)+\varphi(x^{[2]},az_2,\cdots,z_n)\\
&=&a\delta^n\omega(x,az_2,\cdots,z_n)+\theta_E(x)(a)\varphi(x,z_2,\cdots,z_n)+\sum_{i=3}^{n}\theta_E(z_i)(a)\omega(x,z_2,\cdots,\hat{z_i},\cdots,z_n)\\	
&&+\theta_E(x)(a)\varphi\left(z_2,x,z_3,\cdots,\hat{z_i},\cdots,z_n \right)\\	
&&-\sum_{j=3}^n \theta_E(z_j)(a)\omega\left(x,[z_i,z_j],z_2,z_3,\cdots,\hat{z_i},\cdots,\hat{z_j},\cdots,z_n  \right)\\
&=&a\delta^n\omega(x,z_2,\cdots,z_n).
\end{eqnarray*} 
Thus, for any $(\varphi,\omega)\in C_{\rm LR}^{n}(E,M)$, we have  $\dM^{n}_{\rm res}(\varphi,\omega)\in C_{\rm LR}^{n+1}(E,M)$. Furthermore, we have $\dM^{n+1}_{\rm res}\circ \dM^n_{\rm res}=0$. The conclusion follows.

\end{proof}

\subsection{Abelian extensions}
Those extensions have been considered in \cite{Hu} for the characteristic zero case and in \cite{Do2} for the restricted case, using a different approach.

Let $(A,L,\theta_L)$ be a restricted Lie-Rinehart algebra and let
 $(A,M,\theta_M)$ be a strongly abelian restricted Lie-Rinehart algebra (\textit{i.e}, $[m,n]=0~\forall m,n\in M,~\text{and }m^{[p]_M}=0~\forall m\in M$). An \textit{abelian extension} of $(A,L,\theta_L)$ by $(A,M,\theta_M)$ is a short exact sequence of restricted Lie-Rinehart algebras
\begin{equation}\label{eq:ext-seq}
    0\longrightarrow (A,M,\theta_M)\overset{\iota}{\longrightarrow}(A,E,\theta_E)\overset{\pi}{\longrightarrow}(A,L,\theta_L)\longrightarrow 0.
\end{equation}
Since $\iota$ and $\pi$ are morphisms of restricted Lie-Rinehart algebras, it follows that  $\theta_M=0$.

An abelian extension given by a short exact sequence \eqref{eq:ext-seq} is called \textit{split} if there exists an $A$-linear map $\sigma:(A,L,\theta_L)\rightarrow (A,E,\theta_E)$ such that $\pi\circ\sigma=\id$. Such a splitting map $\sigma$ exists if and only if $L$ is projective as $A$-module.

Two abelian extensions of $(A,L,\theta_L)$ by $(A,M,\theta_M)$ are called \textit{equivalent} if there is a restricted Lie-Rinehart algebras morphism $\kappa:(A,E_1,\theta_{1})\rightarrow (A,E_2,\theta_{2})$ such that the following diagram commutes:
\begin{equation}\label{eq:ext-equiv}
\begin{tikzcd}
            &                                                & (A,E_1,\theta_{1}) \arrow[rd, "\pi_1"] \arrow[dd, "\kappa"] &             &   \\
0 \arrow[r] & (A,M,\theta_M) \arrow[ru, "\iota_1"] \arrow[rd, "\iota_2"'] &                                              & (A,L,\theta_L) \arrow[r] & 0. \\
            &                                                & (A,E_2,\theta_{2}) \arrow[ru, "\pi_2"']                     &             &  
\end{tikzcd}
\end{equation}

 \sssbegin{Lemma}
 Consider an abelian extension given by a short exact sequence \eqref{eq:ext-seq}. Then, $(M,\rho)$ is a $(A,L,\theta_L)$-restricted Lie-Rinehart module, the module structure being given by the $A$-linear map $\rho:L\rightarrow\End(M)$ defined by
 \begin{equation}\label{eq:module-ext}
    \rho(x)(m):=\iota^{-1}\bigl([\Tilde{x},\iota(m)]_E\bigl),\quad \forall x\in L,~\forall m\in M,
 \end{equation}
 where $\pi(\Tilde{x})=x.$ 
 \end{Lemma}

 \begin{proof}
Since $M$ is abelian, $\rho$ defines a restricted Lie-module structure, see \cite[Section 5]{EF}. Let $x\in L$, $a\in A$ and $m\in M$. Since $\iota$ is $A$-linear, we have
\begin{align*}
    \rho(x)(am)&=\iota^{-1}\bigl([\Tilde{x},\iota(am)] \bigl)\\
                &=\iota^{-1}\bigl(a[\Tilde{x},\iota(m)]_E+\theta_E(\Tilde{x})(a)\iota(m)\bigl)\\
                &=a\rho(x)(m)+\theta_E(\tilde{x})(a)m.
\end{align*}
Since $\pi$ is a morphism of restricted Lie-Rinehart algebras, we have $$\theta_E(\tilde{x})(a)m=\theta_L\bigl(\pi(\tilde{x})\bigl)(a)m=\theta_L(x)(a)m.$$ Thus, we have $\rho(x)(am)=a\rho(x)(m)+\theta_L(x)(a)m$ and $\rho$ defines a restricted Lie-Rinehart module.
\end{proof}

\sssbegin{Theorem}\label{thmext}
Let $(A,L,\theta_L)$ be a restricted Lie-Rinehart algebra such that $L$ is projective as $A$-module, and let $(A,M,\theta_M)$ be a strongly abelian restricted Lie-Rinehart algebra. Then, the equivalence classes of abelian extensions of $(A,L,\theta_L)$ by $(A,M,\theta_M)$ are classified by $\HH^2_{\rm LR}(L;M)$. 
\end{Theorem}

\begin{proof}
    Consider a split abelian extension given by a short exact sequence \eqref{eq:ext-seq}. The $A$-linear splitting map $\sigma$ defines a restricted Lie-Rinehart module $(M,\rho)$ by taking $\tilde{x}=\sigma(x)$ in Eq. \eqref{eq:module-ext}. Consider the maps $\varphi:L\wedge L\rightarrow M$ and $\w:L\rightarrow M$ defined for all $x,y\in L$ by
    \begin{align}\label{eq:sigmacocycle1}
    \varphi(x,y):=&~\iota^{-1}\Bigl(\bigl[\sigma(x),\sigma(y)\bigl]_E+\sigma\bigl([x,y]_L\bigl)\Bigl);\\\label{eq:sigmacocycle2}
        \w(x):=&~\iota^{-1}\Bigl(\sigma(x)^{[2]_E}+\sigma\bigl(x^{[2]_L}\bigl)\Bigl)
    \end{align}
Let $x,y\in L$ and $a\in A$. Since $\sigma$ is $A$-linear, we have
\begin{align*}
\varphi(x,ay)&=\iota^{-1}\bigl[\sigma(x),\sigma(ay)\bigl]_E+\iota^{-1}\sigma\bigl([x,ay]_L\bigl)\\
    &=\iota^{-1}a\bigl[\sigma(x),\sigma(y)\bigl]_E+\iota^{-1}\theta_E(\sigma(x))(a)\sigma(y)+\iota^{-1}\sigma\bigl(a[x,y]_L+\theta_L(x)(a)y\bigl)\\
    &=a\varphi(x,y)+\iota^{-1}\sigma\bigl(\theta_E(\sigma(x))(a)y+\theta_L(x)(a)y  \bigl)\\
    &=a\varphi(x,y),
\end{align*}
since $\theta_E\circ\sigma=\theta_L\circ\pi\circ\sigma=\theta_L$. Similarly, we have
\begin{align*}
    \w(ax)&=\iota^{-1}\sigma(ax)^{[2]_E}+\iota^{-1}\sigma\bigl((ax)^{[2]_L}\bigl)\\
          &=a^2\w(x)+\iota^{-1}a\theta_E(\sigma(x))(a)\sigma(x)+\iota^{-1}\sigma\bigl(\theta_L(ax)(a)x \bigl)  \\
          &=a^2\w(x).
\end{align*} Thus, we have $(\varphi,\w)\in C_{\rm LR}^2(L;M)$. Moreover, we also have $d^2_{\res}(\varphi,\w)=0$ (see \cite[Lemma 8]{EF}), therefore $(\varphi,\w)\in Z^2_{\LR}(L;M)$.

Conversely, let $(\varphi,\w)\in Z^2_{\LR}(L;M)$. Since the short exact sequence $\eqref{eq:ext-seq}$ splits, we have an isomorphism of $A$-modules $E\cong L\oplus M$. We define the following maps on $L\oplus M$.
\begin{equation*}
 \begin{array}{rllr}
       [x+u,y+v]_{\sigma}&:=&[x,y]_L+\varphi(x,y)+\rho(x)(v)+\rho(y)(u),&\forall x,y\in L,~\forall u,v\in M;\\
       (x+u)^{[2]_{\sigma}}&:=&x^{[2]_L}+\w(x)+\rho(x)(u),&\forall x\in L,~\forall u\in M;\\
       \theta_{\sigma}(x+u)&:=&\theta_L(x),&\forall x\in L,~\forall u\in M.
\end{array}
\end{equation*}
    These maps yield a restricted Lie-Rinehart structure on $L\oplus M$, see \cite[Proposition 3.5]{Do2}.

    Next, we show that the constructions do not depend on the choice of the splitting. Let $\sigma,\sigma':L\rightarrow E$ be two splitting maps and denote by $(\varphi,\w)$ (resp. $(\varphi',\w')$) the corresponding 2-cocycles defined by Eqs. \eqref{eq:sigmacocycle1}, \eqref{eq:sigmacocycle2}. Let $\tau:=\sigma-\sigma'$. Since $\pi\circ\sigma=\pi\circ\sigma'=\id,$ we have that $\tau(x)\in\Ker(\pi),~\forall x\in L$. Thus, $\iota^{-1}\circ\tau(x)\in \iota(M),~\forall x\in L.$
    Since $M$ is abelian, we have for all $x\in L$ and for all $u\in M,$
    \begin{equation}
        \rho(x)(u)=\iota^{-1}\bigl([\sigma(x),\iota(u)]_E\bigl)=\iota^{-1}\bigl([\sigma'(x)+\tau(x),\iota(u)]_E\bigl)=\iota^{-1}\bigl([\sigma'(x),\iota(u)]_E\bigl).
    \end{equation} Therefore, the representation $\rho$ does not depend on the choice of the splitting map. Furthermore, for all $x,y\in L$, we have
    \begin{align*}
        (\varphi+\varphi')(x,y)&=\iota^{-1}\Bigl([\sigma(x),\sigma(y)]_E+\sigma\bigl([x,y]_L\bigl)+[\sigma'(x),\sigma'(y)]_E-\sigma'\bigl([x,y]_L\bigl)\Bigl)\\
        &=\rho(x)\bigl(\iota^{-1}\circ\tau(y)\bigl)+\rho(y)\bigl(\iota^{-1}\circ\tau(x)\bigl)+\iota^{-1}\bigl(\tau([x,y]_L)\bigl)\\
        &=\dM^1_{\rm Rin}(\iota^{-1}\circ\tau)(x,y).
    \end{align*} Moreover, we have
    \begin{align*}
        (\w+\w')(x)&=\iota^{-1}\Bigl(\sigma(x)^{[2]_E}+\sigma(x^{[2]_L})+\sigma'(x)^{[2]_E}+\sigma'(x^{[2]_L})   \Bigl)\\
        &=\iota^{-1}\Bigl(\bigl(\sigma(x)+\sigma'(x)\bigl)^{[2]_E}+[\sigma(x),\sigma'(x)]_E+\tau\bigl(x^{[2]_L}\bigl)
        \Bigl)\\
        &=\underset{=~0}{\underbrace{\bigl(\iota^{-1}\circ\tau(x)\bigl)^{[2]_M}}}+\iota^{-1}\bigl([\sigma(x),\tau(x)]_E\bigl)+\iota^{-1}\circ\tau(x^{[2]_L})\\
        &=\rho(x)\bigl(\iota^{-1}\circ\tau(x)\bigl)+\iota^{-1}\circ\tau(x^{[2]_L})=\delta^1(\iota^{-1}\circ\tau)(x).
    \end{align*} Therefore, we have that $(\varphi+\varphi',\w+\w')=\dM^1_{\rm res}(\iota^{-1}\circ\tau)$. Thus, choosing another splitting leads to cohomologous 2-cocycles.

    It remains to show that equivalent extensions are classified by $\HH^2_{\rm LR}(L;M)$. Let $(A,E_1,\theta_1)$ and $(A,E_2,\theta_2)$ be equivalent extensions of $(A,L,\theta_L)$ given by the diagram \eqref{eq:ext-equiv} and denote by $(\varphi_1,\w_1)$ (resp. $(\varphi_2,\w_2)$) the corresponding 2-cocycles. Since $\pi_2(\kappa\sigma_1+\sigma_2)=0$, we can define an $A$-linear map $\psi:L\rightarrow M,~x\mapsto \iota^{-1}(\kappa\sigma_1+\sigma_2)$. Then, we have $$\dM_{\rm res}^1(\psi)=\bigl(\varphi_1+\varphi_2,\w_1+\w_2\bigl).$$

    Conversely, let $(\varphi_1,\w_1)\in Z^2_{\LR}(L;M)$ and let $(A,E_1,\theta_1)$ be the extension of $(A,L,\theta_L)$ given by $(\varphi_1,\w_1)$.
    For $\psi:L\rightarrow M$ be an $A$-linear map, we consider
    $\varphi_2:=\varphi_1+\dM_{\rm Rin}^1\psi$ and $\w_2:=\w_1+\delta^1\psi.$ We will show that $(A,E_1,\theta_1)$ is equivalent to the extension $(A,E_2,\theta_2)$ given by $(\varphi_2,\w_2)$. We consider the bijective $A$-linear map $\kappa: E_1\rightarrow E_2$ defined by 
    $$\kappa(x+u):=x+u+\psi(x),~\forall x\in L,~\forall u\in M.$$ This map is a restricted morphism of restricted Lie algebras that makes the diagram \eqref{eq:ext-equiv} commute, see \cite[Theorem 5]{EF}. Moreover, it also satisfies
    $$\theta_2\bigl(\kappa(x+u)\bigl)(a)=\theta_1(x+u)(a),~\forall x\in L,~\forall a\in A,~\forall u\in M.$$ Therefore, $\kappa$ is an isomorphism of restricted Lie-Rinehart algebras and the two extensions $(A,E_1,\theta_1)$ and $(A,E_2,\theta_2)$ are equivalent.
\end{proof}

\subsection{A comparison with Dokas' work \cite{Do2}}


In \cite{Pa}, Pareigis introduced the cohomology of restricted Lie algebras, which differs from the Hochschild cohomology introduced in 1954. In degrees 1 and 2, both cohomologies differ, but coincide in degrees higher than 2. In fact $\mathrm{H}^2_{\mathrm{Pa}}(L;M)$  classifies extensions of the restricted Lie algebra  $L$ by an abelian restricted module $M$ admitting a non-necessarily zero $p$-map and on which $L$ acts upon by $p$-derivations, while Hochschild (and Evans-Fuchs) $\mathrm{H}^2_{\mathrm{res}}(L;M)$  classifies central extensions of $L$ by a {\it strongly abelian} restricted module $M$. 
  In \cite{Do1}, Dokas recovered Pareigis constructions using Quillen-Barr-Beck methods, which are generally applicable to algebraic structures, especially restricted Lie Rinehart, see \cite{Do2}. Therefore, Dokas constructions classify extensions of restricted Lie-Rinehart algebras by abelian non-trivial modules with $p$-maps, while ours (see Theorem \ref{thmext} ) classifies extensions by strongly abelian modules.\\
  
 {\bf Open problems:} (i) It would be interesting to build a dictionary that compares the two approaches. It remains a gap in the literature to find a correspondence between Hochschild-Evans-Fuchs cohomology in degrees 1 and 2 and that constructed using the Quillen-Barr-Beck method.\\
(ii) Compare the Poisson and restricted Lie Rinehart cohomology introduced in this paper with that based on the Quillen-Barr-Beck method.

\black

\section{Cohomology and deformations of restricted Poisson algebras}\label{sec:poisson}

An associative commutative $\K$-algebra $(A,\cdot)$ (not necessarily unital) is called \textit{Poisson algebra} if it is equipped with a bilinear map $\{-,-\}:A\times A\rightarrow A$ such that $\bigl(A,\{-,-\}\bigl)$ is a Lie algebra and moreover, we have
\begin{equation}\label{eq:poisson-definition}
    \{a\cdot b,c\}=a\cdot \{b,c\}+b\cdot \{a,c\},\quad \forall a,b,c\in A.
\end{equation}
When no confusion is possible, the associative commutative product $\cdot$ will be denoted by juxtaposition. We denote such a Poisson algebra by a triple $\bigl(A,\cdot,\{-,-\}\bigl)$.

\subsection{Restricted Poisson algebras for $p=2$}\label{def:restrictedpoisson} A Poisson algebra $\bigl(A,\cdot,\{-,-\}\bigl)$ in characteristic $2$ is called a\textit{ weakly restricted Poisson algebra} if there exists a map $(-)^{\{2\}}:A\rightarrow A$ such that $\bigl(A,\{-,-\},(-)^{\{2\}}\bigl)$ is a restricted Lie algebra. Furthermore, following \cite{PS}, a weakly restricted Poisson algebra $A$ is called a \textit{restricted Poisson algebra} if we have
  \begin{equation}\label{eq:resPA-2map}
(xy)^{\{2\}}=x^2y^{\{2\}}+ y^2x^{\{2\}}+xy\{x,y\},~\forall x,y\in A.
  \end{equation}
\sssbegin{Lemma}
    Let $\bigl(A,\cdot,\{-,-\},(-)^{\{2\}} \bigl)$ be a weakly restricted Poisson algebra. Then, we have
    \begin{equation}
        \ad_{xy^{\{2\}}}=\ad_{x^2y^{\{2\}}}+\ad_{y^2x^{\{2\}}}+\ad_{xy\{x,y\}},\quad \forall x,y\in A.
    \end{equation}
\end{Lemma}

\begin{proof} Straightforward computations.
\end{proof}

\sssbegin{Remark}
 In \cite{BK2}, the analog of Condition \eqref{eq:resPA-2map} for $p\geq 3$ reads
 \begin{equation}\label{eq:byz}
    (xy)^{\{p\}}=x^2y^{\{p\}}+ y^2x^{\{2\}}+\Phi_p(x,y),~\forall x,y\in A,
 \end{equation} where
 \begin{equation*}
        \Phi_p(x,y)=(x^p+y^p)\sum_{1\leq i\leq p-1}s_i(x,y)-\frac{1}{2}\sum_{1\leq i\leq p-1}s_i(x^2,y^2)+\sum_{1\leq i\leq p-1}s_i(x^2+y^2, 2xy),\quad \forall x,y\in A.
 \end{equation*}
 In the case where $p\geq 3$, one can show that Equation \eqref{eq:byz} is equivalent to
 \begin{equation}\label{eq:byz2}
    (x^2)^{\{p\}}=2x^px^{\{p\}},\quad\forall x\in A,
 \end{equation} which is the condition that appears in \cite[Definition 3.4]{BYZ}. In characteristic 2, Equation \eqref{eq:byz2} reduces to $(x^2)^{\{2\}}=0$, which  is not equivalent to Equation \eqref{eq:resPA-2map}. Hence, the definition given in Equation \eqref{eq:resPA-2map}.
\end{Remark}

\sssbegin{Lemma} \textup{(}see also \cite[Theorem 4.7]{BYZ}\textup{)} Let $\bigl(A,\cdot,\{-,-\},(-)^{\{2\}} \bigl)$ be a weakly restricted Poisson algebra, and let $\mathcal{B}$ be a linear basis of $A$. Then, $\bigl(A,\cdot,\{-,-\},(-)^{\{2\}} \bigl)$ is a restricted Poisson algebra if and only if Equation \eqref{eq:resPA-2map} holds for every pair of elements of $\mathcal{B}$.
\end{Lemma}

\begin{proof}
    Suppose that
    \begin{equation}
        (ab)^{\{2\}}=a^2b^{\{2\}}+ b^2a^{\{2\}}+ab\{a,b\},~\forall a,b\in \mathcal{B}.
    \end{equation} Let $a\in \mathcal{B}$ and consider the set
    \begin{equation*}
        R_a=\bigl\{x\in A,~ (ax)^{\{2\}}=a^2x^{\{2\}}+ x^2a^{\{2\}}+ax\{a,x\} \bigl\}.
    \end{equation*} Then, $R_a$ is a linear subspace of $A$ containing the basis $\mathcal{B}$. Therefore, $R_a=A,~\forall a\in \mathcal{B}$. Similarly, we show that $L_a=A,~\forall a\in\mathcal{B},$ where
    $$L_a=\bigl\{x\in A,~ (xa)^{\{2\}}=x^2a^{\{2\}}+ a^2x^{\{2\}}+xa\{x,a\} \bigl\}.$$ Hence the conclusion.
\end{proof}

\sssbegin{Remark} Let $\bigl(A,\cdot,\{-,-\},(-)^{\{2\}} \bigl)$ be a restricted Poisson algebra and denote by $Z(A)$ its center (as a Lie algebra). Following \cite{BYZ}, a map $\gamma:A\rightarrow Z(A)$ is called a \textit{central Frobenius derivation} if
\begin{equation}
    \gamma(\lambda x+y)=\lambda^2\gamma(x)+\gamma(y);\quad\text{and}\quad\gamma(xy)=x^2\gamma(y)+y^2\gamma(x),\quad \forall x,y\in A.
\end{equation} Then, a map $(-)^{\{2\}}+\gamma$ defines another restricted Poisson structure on $\bigl(A,\cdot,\{-,-\}\bigl)$ if and only if $\gamma$ is a central Frobenius derivation of $A$.
\end{Remark}

A linear map $\phi:\bigl(A,\cdot,\{-,-\}_A,(-)^{\{2\}_A}\bigl)\rightarrow \bigl(B,\cdot,\{-,-\}_B,(-)^{\{2\}_B}\bigl) $ is a {\it restricted Poisson morphism} if for all $x,y\in A$,
\begin{align}
    \phi(\{x,y\}_A)&=\{\phi(x),\phi(y)\}_B,\\
    \phi(x^{\{2\}_A})&=\phi(x)^{\{2\}_B},\\
    \phi(x\cdot y)&=\phi(x) \phi(y).
\end{align}

A linear map $D:\bigl(A,\cdot,\{-,-\}_A,(-)^{\{2\}_A}\bigl)\rightarrow \bigl(A,\cdot,\{-,-\}_A,(-)^{\{2\}_A}\bigl)$ is a {\it restricted Poisson derivation} if for all $x,y\in A$,
\begin{align}
    D(\{x,y\})&=\{D(x),y\}+\{x,D(y)\},\\
    D(x^{\{2\}})&=\{x,D(x)\},\\
    D(x\cdot y)&=D(x)\cdot y+x\cdot D(y).
\end{align}

Recall that a {\it Hamiltonian derivation} of $A$ is a restricted Poisson derivation $D$ induced by an element $a\in A$ such that 
$$D(b):=\{a,b\},\quad \forall~b\in A.$$

We denote the set of Poisson derivations and Hamiltonian derivations of the restricted Poisson algebra $A$ by {\rm PDer($A$)} and {\rm HDer($A$)}, respectively.

The following Proposition is a restricted version of \cite[Proposition 13.3.8]{LV}.
\sssbegin{Proposition}
    Let $(A,L,\theta)$ be a restricted Lie-Rinehart algebra. Then, the tuple $\bigl(A\oplus L,\cdot,\{-,-\},(-)^{\{2\}}\bigl)$ is a restricted Poisson algebra, where
\begin{equation}
\begin{array}{rlll}(a+x)\cdot(b+y)&:=&ab+ay+bx;&~\forall a,b\in A,~\forall x,y\in L;\\ [2mm]
\{a+x,b+y\}&:=&[x,y]_L+\theta(x)(b)+\theta(y)(a),&~\forall a,b\in A,~\forall x,y\in L;  \\[2mm]
(a+x)^{\{2\}}&:=&x^{[2]_L}+\theta(x)(a),&~\forall a\in A,~\forall x\in L.
\end{array}
\end{equation}
\end{Proposition}
\begin{proof}
    Let $a,b\in A $ and $x,y\in L$. We have
    \begin{align*}
        \{a^{\{2\}},x\}+\{a,\{a,x\}\}&=\{a,\theta(x)(a)\}=0;\quad\text{and}\\
         \{x^{\{2\}},a\}+\{x,\{x,a\}\}&=\theta(x^{[2]})(a)+\theta(x)^2(a)=0. \end{align*}
         Thus, $(A\oplus L, \{-,-\},(-)^{\{2\}})$ is a restricted Lie algebra. Moreover, we have
         \begin{align*}
            (ax)^{\{2\}}+a^2x^{\{2\}}+ax\cdot\{a,x\}&=2a^2x^{[2]_L}+\theta(ax)(a)x+(\theta(x)(a)a)x\\
            &=\theta(ax)(a)x+(a\theta(x)(a))x \\
            &=0.          
         \end{align*}The remaining identities are either trivial or covered by \cite[Proposition 13.3.8]{LV}. Therefore, $\bigl(A\oplus L,\cdot,\{-,-\},(-)^{\{2\}}\bigl)$ is a restricted Poisson algebra.
\end{proof}

\subsubsection{Examples from deformation quantization}\label{sec:exdefoquant}

Let $A$ be an associative algebra over an arbitrary field $\K$, and let $M$ be an $A$-bimodule. The objective of this subsection is to provide a generic method of constructing restricted Poisson algebras based on formal deformations of associative algebras. For completeness, we recall Hochschild's cohomology of associative algebras, see \cite{H}. The cochains space are defined by
\begin{equation}
\begin{array}{rlll}
C_{\mathrm H}^n(A;M)&:=&\Hom(A^{\otimes n},M)&\forall n>0;\\ [2mm]
C_{\mathrm H}^0(A;M)&:=& M  \\[2mm]
C_{\mathrm H}^n(A;M)&:=&0~&\forall n<0.
\end{array}
\end{equation}
The differential maps $\dM_{\mathrm H}: C^n_{\mathrm H}(A;M)\rightarrow C^{n+1}_{\mathrm H}(A;M)$  are defined by (for all $a_1,\ldots,a_{n+1}\in A$)
\[
\begin{array}{lcl}(\dM_{\rm H}^n\varphi)(a_1,a_2,...,a_{n+1})&:=&\displaystyle a_1\varphi(a_2,...,a_{n+1})+\sum_{i=1}^{n}(-1)^i\varphi(a_1,...,a_ia_{i+1},...,a_{n+1})\\[2mm]
&&+(-1)^{n+1}\varphi(a_1,...,a_n)a_{n+1}.
\end{array}
\]

We have $d_{\mathrm H}^{n+1}\circ d_{\mathrm H}^n=0$. We denote by $Z_{\mathrm H}^n(A;M)=\text{Ker}(d_{\mathrm H}^n)$ the Hochschild $n$-cocycles and by $B_{\mathrm H}^n(A;M)=\text{Im}(d_{\mathrm H}^{n-1})$ the Hochschild $n$-coboundaries.
The Hochschild cohomology groups are defined by
\[\HH_{\mathrm H}^n(A;M):=Z_{\mathrm H}^n(A;M)/B_{\mathrm H}^n(A;M).\]

\noindent\textbf{Notations.} Let $A$ be an associative algebra. In this subsection, we denote by $A^-$ the restricted Lie algebra which is given by the commutator bracket and the $2$-map given by $a^{[2]}:=a^2$, for all $a\in A$. Moreover, for any bilinear map $\mu:A\times A\rightarrow M$, we denote by $\mu^-$ the commutator of $\mu$ defined by
$$\mu^-(a,b):=\mu(a,b)-\mu(b,a),~\forall a,b\in A;$$ and by $\w_{\mu}$ the map 
$$\omega_{\mu}(a):=\mu(a,a),~\forall a\in A.$$


\sssbegin{Lemma}\label{leibniz}
    Let $A$ be an associative commutative algebra and let $\mu\in Z_{\mathrm H}^2(A;A)$.  Then, $\mu^-\in\Der(A)$.    
\end{Lemma}

\begin{proof}
    Let $a,b,c\in A$. Using the $2$-cocycle condition, we have 
    \begin{align*}
        \mu^-(ab,c)&=\mu(ab,c)+\mu(c,ab)\\
                   &=a\mu(b,c)+\mu(a,bc)+\underline{\mu(a,b)c+c\mu(a,b)}+\mu(ca,b)+\mu(c,a)b.
    \end{align*}
    As $A$ is commutative, the underlined terms vanish. Using the $2$-cocycle condition again, we get
$$\mu(a,bc)+\mu(ca,b)=\mu(a,cb)+\mu(ac,b)=a\mu(c,b)+\mu(a,c)b.$$  Therefore, we finally obtain
    \begin{align*}
         \mu^-(ab,c)&=\mu(ab,c)+\mu(c,ab)\\
                    &=a\mu(b,c)+a\mu(c,b)+\mu(a,c)b+\mu(c,a)b\\
                    &=a\mu^-(b,c)+\mu^-(a,c)b.
    \end{align*} Thus, $\mu^-$ is a derivation of the associative algebra $A$.
\end{proof}

A \textit{formal deformation} of an associative algebra $A$ is an associative product $\mu_t$ on the formal power series algebra $A_t:=A[[t]]$ given by
 \begin{equation}
    \mu_t(a,b)=ab+\sum_{i\geq 1}t^i\mu_i(a,b),~\forall a,b\in A,
\end{equation}
where the maps $\mu_i:A\times A\rightarrow A$ are bilinear. The associativity of $\mu_t$ is equivalent to 
 \begin{equation}\label{eqdeformation}
    \sum_{i=0}^k\mu_i\bigl(a,\mu_{k-i}(b,c)\bigl)=\sum_{i=0}^k\mu_i\bigl(\mu_{k-i}(a,b),c\bigl)~\forall a,b,c\in A,~\forall k\geq 0.
\end{equation}
In particular, $\mu_1$ is a $2$-cocycle of the Hochschild cohomology, see \cite{Ge}. Formal deformations of restricted Poisson algebras are investigated in Section \ref{sec:defopoisson}.

In the sequel, let $A$ be an associative commutative algebra and let $(A_t,\mu_t)$ be a formal deformation of $A$.

\sssbegin{Lemma} For every  $a,b\in A$, define $\{a,b\}:=\mu_1^-(a,b).$ Then, $\{-,-\}$ is a Poisson bracket on $A$.
\end{Lemma}

\begin{proof}
   We show that $\{-,-\}$ satisfies the Jacobi identity.
    Let $a,b,c\in A$ and consider the map $\obs_2\in C_{\mathrm H}^3(A;A)$ defined by 
   $$\obs_2(a,b,c):=\mu_1(a,\mu_1(b,c))+\mu_1(\mu_1(a,b),c),~\forall a,b,c\in A.$$ This map represents the obstruction to the integrability of $\mu_2$. In particular, if $\mu_t$ is associative, then $\obs_2=\dM^2_{\rm H}\mu_2$, see \cite{Ge}. We have

   \begin{align*}
    &\{a,\{b,c\}\}+\{b,\{c,a\}\}+\{c,\{a,b\}\}\\
    =~&\mu_1(a,\mu_1(b,c))+\mu_1(a,\mu_1(c,b))+\mu_1(\mu_1(b,c),a)+\mu_1(\mu_1(c,b),a)\\
    &~~+\mu_1(b,\mu_1(c,a))+\mu_1(b,\mu_1(a,c))+\mu_1(\mu_1(c,a),b)+\mu_1(\mu_1(a,c),b)\\
    &~~+\mu_1(c,\mu_1(a,b))+\mu_1(c,\mu_1(b,a))+\mu_1(\mu_1(a,b),c)+\mu_1(\mu_1(b,a),c)\\
    =~&\obs_2(a,b,c)+\obs_2(a,c,b)+\obs_2(b,c,a)+\obs_2(c,b,a)+\obs_2(c,a,b)+\obs_2(b,a,c)\\
    =~&\dM^2_{\mathrm H}\mu_2(a,b,c)+\dM^2_{\mathrm H}\mu_2(a,c,b)+\dM^2_{\mathrm H}\mu_2(b,c,a)+\dM^2_{\mathrm H}\mu_2(c,b,a)+\dM^2_{\mathrm H}\mu_2(c,a,b)+\dM^2_{\mathrm H}\mu_2(b,a,c)\\
    =~&0.
   \end{align*} The last equality is obtained by expanding each coboundary and canceling terms two by two in the expansion. 

    Moreover, $\{-,-\}$ yields a derivation of the associative product of $A$ by Lemma \ref{leibniz}. Thus, it is a Poisson bracket.
\end{proof}
\sssbegin{Proposition}\label{propdef}
Let $A$ be an associative commutative algebra, and let $(A_t,\mu_t)$ be a formal deformation of $A$ such that $\mu_1(a^2,b)=0 \text{ and }\mu_2^-(a^2,b)=0~\forall a,b\in A$. Then, $(A,\mu_1^-,\omega_{\mu_1})$ is a restricted Poisson algebra.
\end{Proposition}

\begin{proof}
  Since $(A_t,\mu_t)$ is an associative algebra, it follows that $(A_t,\mu_t^-,\w_{\mu_t})$ is a restricted Lie algebra. Therefore, we have
  \begin{align*}
    \mu_t^-\bigl(a,\mu_t^-(a,b)\bigl)+\mu_t^-\bigl(\mu_t(a,a),b\bigl)=0,\quad \forall a,b\in A.
  \end{align*}
  By expanding this relation and collecting the coefficients of $t^2$, we obtain
  $$ \mu_1^-\bigl(a,\mu_1^-(a,b)\bigl)+\mu_2^-(a^2,b)+\mu_1^-\bigl(\mu_1(a,a),b\bigl),\quad \forall a,b\in A.$$ Since $\mu_2^-(a^2,b)=0$, we obtain
  $$ \mu_1^-\bigl(a,\mu_1^-(a,b)\bigl)=\mu_1^-\bigl(\w_{\mu_1}(a),b\bigl),~\forall a,b\in A. $$ Thus, $(A,\mu_1^-,\omega_{\mu_1})$ is a restricted Lie algebra. It remains to show the identity \eqref{eq:resPA-2map}. Using the Hochschild 2-cocycle condition, we have
  \begin{align*}
    \w_{\mu_1}(ab)=\mu_1(ab,ab)&=a\mu_1(b,ab)+\mu_1(a,ab^2)+\mu_1(a,b)b^2\\
    &=a\mu_1(a,b^2)+\mu_1(a^2,b^2)+\mu_1(a,a)b^2\\&~+ab\mu_1(b,a)+a\mu_1(b^2,a)+a^2\mu_1(b,b)+\mu_1(a,b)ab\\
    &=\mu_1(a,b)b^2+ab\mu_1(b,a)+a^2\mu_1(b,b)+\mu_1(a,b)ab\\
    &=a^2\w_{\mu_1}(b)+b^2\w_{\mu_1}(a)+ab\mu_1^-(a,b).
  \end{align*}
  Note that since $\mu_1^-$ is a derivation, we have $\mu_1^-(a^2,b)=0,~\forall a,b\in A.$ Then, we have $\mu_1(b,a^2)=\mu_1(a^2,b)=0,~\forall a,b\in A$ in the above computation.
  Therefore, $(A,\mu_1^-,\omega_{\mu_1})$ is a restricted Poisson algebra.
\end{proof}

\subsection{Cohomology of restricted Poisson algebra for $p=2$} In this Section, we introduce a cohomology theory for restricted Poisson algebras for $p=2$. This cohomology is constructed using the restricted cohomology of restricted Lie algebras described in Section \ref{sec:rescoho2}.
\subsubsection{Skew-symmetric $k$-derivations} Let $A$ be a commutative associative algebra in characteristic $2$ and let $k\geq0$. A skew-symmetric $k$-linear map $\varphi\in \Hom(\wedge^kA,A)$ is called a \textit{skew-symmetric $k$-derivation of $A$} if $\varphi$ is a derivation in each of its arguments. Due to skew-symmetry, this is equivalent to requiring that $\varphi$ is a derivation in its first argument, that is, 
\begin{equation}\label{eq:derivation}
   \varphi(xy,z_2,\cdots,z_k)=x\varphi(y,z_2,\cdots,z_k)+y \varphi(x,z_2,\cdots,z_k),~\forall x,y,z_2,\cdots,z_k\in A.
\end{equation}
The vector space of all skew-symmetric $k$-derivations of $A$ is denoted by $\mathfrak{X}^k(A)$ and we introduce the graded vector space 
$$\mathfrak{X}^\bullet(A):=\bigoplus_{k\geq 0}\mathfrak{X}^k(A).$$
\subsubsection{Restricted Poisson cochains} Let $(A,\cdot,\{-,-\},(-)^{\{2\}})$ be a restricted Poisson algebra in characteristic $2$. For $n=0,1$, we set $C^0_{\rm PA}(A)=A$ and $C^1_{\rm PA}(A)=\frak{X}^1(A)$. For $n\geq 2$, a pair  $(\varphi,\omega)\in C_{\res}^n(A;A)$ is a restricted Poisson algebra $n$-cochain if for all $\lambda\in \K$ and all $x,y,z_2,\cdots,z_{n-1}\in A$, we have (where  $2\leq i\leq n-1$):
\begin{align}
\label{eq:res PA1}		
\omega( x y, z_2,\cdots,z_{n-1})=&~x^2\omega(y,z_2,\cdots,z_{n-1})+y^2\omega(x,z_2,\cdots,z_{n-1})\\
\nonumber&+xy\varphi(x,y,z_2,\cdots,z_{n-1}),\\
\label{eq:res PA2}			
\omega( x, z_2,\cdots,z_i z'_i,\cdots,z_{n-1})=&~z_i \omega(x,z_2,\cdots,z'_i,\cdots,z_{n-1})\\\nonumber&~+z'_i\omega(x,z_2,\cdots,z_i,\cdots,z_{n-1}).
\end{align}
The space of all pairs $(\varphi,\w)\in C_{\res}^n(A;A)$ satisfying \eqref{eq:res PA1} and \eqref{eq:res PA2} is denoted by $C^n_{\rm PA}(A)$. The coboundary operators
${\rm d}_{\rm PA}^n:C^n_{\rm PA}(A)\rightarrow C^{n+1}_{\rm PA}(A)$ for
$n\geq 0$ are induced by the restricted coboundary operators for restricted Lie algebras, see Section \ref{sec:diff.op.p=2}, where $A$ is seen as a module over itself with the adjoint action.

We denote by
$Z_{\rm PA}^n(A)=\Ker(\dM^n_{\rm PA})$ the restricted Poisson $n$-cocycles and $B_{\rm PA}^n(A)=\text{Im}(\dM^{n-1}_{\rm PA})$ the restricted Poisson $n$-coboundaries. 

\sssbegin{Proposition}\label{prop:cohopoisson}
Let $(A,\cdot,\{-,-\},(-)^{\{2\}})$ be a restricted Poisson algebra in characteristic $2$, and let $n\geq 0$. For all $(\varphi,\omega)\in  C^{n}_{\rm PA}(A)$, we have  $\dM^{n}_{\rm PA}(\varphi,\omega)\in  C^{n+1}_{\rm PA}(A)$. 
\end{Proposition}
Thus, we obtain a well-defined cochain complex $(\bigoplus_{n\geq0} C^{n}_{\rm PA}(A),\dM^n_{\rm PA})$ and the associated $n$-th cohomology group $\HH_{\rm PA}^n(A)$ is well-defined.

\begin{proof}
For all $x,y,z_2,\cdots,z_n\in A$, we have
\[
\begin{array}{l}
  \delta^n\omega(xy,z_2,\cdots,z_n)\\[1mm]
  =\{xy,\varphi(xy,z_2,\cdots,z_n)\}+ \displaystyle\sum_{i=2}^{n}\{z_i,\omega(xy,z_2,\cdots,\hat{z_i},\cdots,z_n)\} +\varphi((xy)^{\{2\}},z_2,\cdots,z_n)\\[2mm]
  +\displaystyle \sum_{i=2}^{n}\varphi\left(\{xy,z_i\},xy,z_2,\cdots,\hat{z_i},\cdots,z_n \right) +\displaystyle \sum_{2\leq i<j\leq n}\omega\left(xy,\{z_i,z_j\},z_2,\cdots,\hat{z_i},\cdots,\hat{z_j},\cdots,z_n  \right)\\[5mm]
=x^2\delta^n\omega(y,z_2,\cdots,z_n)+y^2\delta^n\omega(x,z_2,\cdots,z_n)+xy\,d_{\ce}^n \varphi(x,y,z_2,\cdots,z_n)\\[2mm]
\quad +2x\{x,y\}\, \varphi(y,z_2,\cdots,z_n)+2y\{x,y\}\, \varphi(x,z_2,\cdots,z_n)+2x^{\{2\}}y\, \varphi(y,z_2,\cdots,z_n)\\[2mm]
\quad +2y^{\{2\}}x\, \varphi(x,z_2,\cdots,z_n)+\displaystyle 2\sum_{i=2}^n\Big(\{z_i,x\}x\omega(y,z_2,\cdots,\hat{z_i},\cdots,z_n)\\[2mm]
\quad +\{z_i,y\}y\, \omega(x,z_2,\cdots,\hat{z_i},\cdots,z_n) + \{z_i,x\}y\, \varphi(x,y,z_2,\cdots,z_n)+\{z_i,y\}x\, \varphi(x,y,z_2,\cdots,z_n)    \Big)\\[2mm]
=x^2\delta^n\omega(y,z_2,\cdots,z_n)+y^2\delta^n\omega(x,z_2,\cdots,z_n)+xy\, d_{\ce }^n \varphi(x,y,z_2,\cdots,z_n).
\end{array}
\]
Furthermore, since $\delta^n\omega$ is skew-symmetric in the $z's$, we only need to check that 
$$\delta^n\omega(x,z_2 z'_2,\cdots,z_n)=z_2\delta^n\omega(x,z'_2,\cdots,z_n)+z'_2 \delta^n\omega(x,z_2,\cdots,z_n).$$
Indeed, 
\[
\begin{array}{l}
\delta^n\omega(x,z_2 z'_2,\cdots,z_n)\\[2mm]
=\{x,\varphi(x,z_2 z'_2,\cdots,z_n)\}+\{z_2z_2',\omega(x,z_3,\cdots,z_n)\}\\[2mm] \quad + \displaystyle \sum_{i=3}^{n}\{z_i,\omega(x,z_2z_2',\cdots,\hat{z_i},\cdots,z_n)\}+\varphi(x^{\{2\}},z_2z_2',\cdots,z_n)+\varphi\left(\{x,z_2z_2'\},x,z_3,\cdots,z_n \right)\\[2mm]
\quad +\displaystyle \sum_{i=3}^{n}\varphi\left(\{x,z_i\},x,z_2z_2',\cdots,\hat{z_i},\cdots,z_n \right)
+\displaystyle 
 \sum_{j=3}^n\omega\left(x,\{z_2z_2',z_j\},z_3,\cdots,\hat{z_j},\cdots,z_n  \right)\\[2mm]
\quad + \displaystyle \sum_{3\leq i<j\leq n}\omega\left(x,\{z_i,z_j\},z_2z_2',\cdots,\hat{z_i},\cdots,\hat{z_j},\cdots,z_n  \right)\\[3mm]
=z_2\delta^n\omega(x, z'_2,\cdots,z_n)+z'_2\delta^n\omega(x, z_2,\cdots,z_n)+2\{x,z_2\}\varphi(x,z_2',\cdots,z_n)+ 2\{x,z'_2\}\varphi(x,z_2,\cdots,z_n)\\[2mm]
\quad +\displaystyle 2\sum_{i=3}^n\big( \{z_i,z'_2\}\omega(x,z_2,\cdots,\hat{z_i},\cdots,z_n)+ \{z_i,z_2\}\omega(x,z'_2,\cdots,\hat{z_i},\cdots,z_n)   \big)\\[2mm]
=z_2\delta^n\omega(x, z'_2,\cdots,z_n)+z'_2\delta^n\omega(x, z_2,\cdots,z_n).
    \end{array}
    \]
As $\dM_{\rm PA}^n$ is induced by the coboundary operator for the restricted Lie algebra $(A,\{-,-\},(-)^{\{2\}})$, we have $\dM^{n}_{\rm PA}(\varphi,\omega)\in  C^{n+1}_{\rm PA}(A)$ and $\dM^{n+1}_{\rm PA}\circ \dM^n_{\rm PA}=0$.
\end{proof}
The cohomology groups of order $0$ and $1$ are given by $$\HH_{\rm PA}^0(A)=\{a\in A\mid \{a,b\}=0,\quad\forall~b\in A\}$$
and 
$$\HH_{\rm PA}^1(A)=\frac{{\rm PDer}(A)}{{\rm HDer}(A)}.$$

Describing the cohomology complex of restricted Poisson algebras in characteristic  $p$  for general $ p > 0$, not just in the case $ p = 2 $ seems to be a difficult problem. In this sense, our work provides a first step toward building such a complex in arbitrary positive characteristic.
\subsection{Comparison between the two cohomologies}
In this section, all commutative associative $\K$-algebras $A$ are supposed to have a unit denoted by $1_A$. In particular, if $(A,\cdot,\{-,-\},(-)^{\{2\}})$ is a restricted Poisson algebra, we have $\{1_A,-\}=1_A^{\{2\}}=0$. 

\subsubsection{K\"{a}hler differentials} We recall the construction of K\"{a}hler differentials, see \cite[Section 3.2.1]{LPV} for more details. Let $A$ be an \textit{unital} associative commutative $\K$-algebra. The module of K\"{a}hler differentials of $A$, denoted by $\Omega^1(A)$, is the free $A$-module generated by elements of the form $\dM x$ for $x\in A$ modulo the submodule generated by elements of the form
\begin{equation}
    \dM (xy)+y\dM x+x\dM y;\quad \dM(x+y)+\dM x+\dM y;\quad \text{and }~ \dM\lambda, \quad \forall x,y\in A,~\forall \lambda\in\K.
\end{equation}
One can consider $\dM$ as a map $\dM:A\rightarrow \Omega^1(A)$. In that case, the pair $(\Omega^1(A),\dM)$ satisfies the following universal property:  for every derivation $D:A\rightarrow A$, there exists a unique $A$-linear map $\hat{D}:\Omega^1(A)\rightarrow A$ such that $\hat{D}\circ d=D$, that is, 
\begin{equation}\label{eq:diag-uni}
\xymatrix{
A \ar[rr]^(.4){\dM}\ar[d]_{D} && \Omega^1(A)
\ar@{.>}[dll]^(.5){\exists !~\hat{D}}\\
A},
\end{equation}
where $A$ is seen as a module over itself using its multiplication. By this universal property, we have a natural isomorphism of $A$-modules
\begin{equation}\label{eq:isomorphism deri}
\Der(A)\cong\Hom_{A}(\Omega^1(A),A).
\end{equation}
Denote by $m:A\otimes A\rightarrow A$ the multiplication of $A$ and by $I=\Ker(m)$. Consider the $A$-modules map
\begin{equation}
    \varphi:I/I^2\rightarrow \Omega^1(A),\quad x\otimes1-1\otimes x\mapsto \dM x.
\end{equation} Then, $\varphi$ induces an isomorphism of $A$-modules between $I/I^2$ and $\Omega^1(A)$, see e.g. \cite[Section 3]{Hu}. Actually,  using the fact that 
\begin{equation}
  x\otimes y+y\otimes x\equiv xy\otimes 1_A+1_A\otimes xy \quad \text{mod}~I^2,\quad\forall x,y\in I,
\end{equation} one can show that $\dM$ (viewed as a map $A\rightarrow \Omega^1(A)$) satisfies the Leibniz identity if and only if $\varphi$ is a morphism of $A$-modules.

The following construction is due to Huebschmann, see \cite[Lemma 3.5 and Theorem 3.8]{Hu}. Let $(A,\cdot,\{-,-\})$ be an ordinary Poisson algebra and $\Omega^1(A)$ its module of K\"{a}hler differentials. The map
\begin{align}
    \pi:\Omega^1(A)\otimes_A\Omega^1(A)\rightarrow A,\quad x\dM u\otimes y\dM v\mapsto xy\{u,v\}
\end{align} is $A$-linear and symmetric. Moreover, the map $\pi$ induces a map
\begin{align}
    \pi^{\sharp}:\Omega^1(A)\rightarrow \Hom_A\bigl(\Omega^1(A),A\bigl)\cong \Der(A),\quad x\dM u\mapsto x\{u,-\}.
\end{align} 
The module of K\"{a}hler differentials $\Omega^1(A)$ is a Lie algebra with the bracket
\begin{equation}
    [x\dM u,y\dM v]_{\Omega^1(A)}:=x\{u,y\}\dM v+y\{x,v\}\dM u+xy\dM\{u,v\},\quad \forall x,y,u,v\in A.
\end{equation}
Then, $\bigl(A,\Omega^1(A),\pi^{\sharp} \bigl)$ is a Lie-Rinehart algebra and $\pi^{\sharp}$ is a morphism of Lie-Rinehart algebras, see \cite[Theorem 3.8]{Hu}. The next result is a restricted version for $p=2$ of \cite[Theorem 3.8]{Hu}. The case of $p\geq 3$ was done in \cite[Theorem 8.2]{BYZ}.

\sssbegin{Theorem}\label{thm:RLRfromPoisson}
    Let $(A,\cdot,\{-,-\},(-)^{\{2\}})$ be a restricted Poisson algebra in characteristic 2, and suppose that the module of K\"{a}hler differentials is free as an $A$-module. Define a map $(-)^{[2]}:\Omega^1(A)\rightarrow \Omega^1(A)$ by
    \begin{equation}
        (x\dM u)^{[2]}=x^2\dM(u^{\{2\}})+x\{u,x\}\dM u,\quad \forall~ x\dM u\in \Omega^1(A).
    \end{equation} Then, $\bigl(\Omega^1(A),[-,-]_{\Omega^1(A)},(-)^{[2]}\bigl)$ is a restricted Lie algebra and $\bigl(A, \Omega^1(A),\pi^{\sharp}\bigl)$ is a restricted Lie-Rinehart algebra.
\end{Theorem}
\begin{proof}
First, we show that the $2$-map is well-defined on $\Omega^1(A)$. Let $x,u,v\in A$. We have
\begin{align*}
    \bigl(x\dM(uv)\bigl)^{[2]}=&~x^2\dM\bigl((uv)^{\{2\}}\bigl)+x\{uv,x\}\dM(uv)\\
    =&~x^2\dM(u^2v^{\{2\}})+x^2\dM(v^2u^{\{2\}})+x^2uv\dM\{u,v\}+x^2u\{u,v\}\dM v+x^2v\{u,v\}\dM u\\
    &~+xu^2\{v,x\}\dM v+xuv\{u,x\}\dM v+xvu\{v,x\}\dM u+xv^2\{u,x\}\dM u\\
    =&~(xu\dM v)^{[2]}+(xv\dM u)^{[2]}+[xu\dM v,xv\dM u]_{\Omega^1(A)}\\
    =&~(xu\dM v+xv\dM u)^{[2]}.
\end{align*} 

Additionally, we have $\bigl(x\dM(u+v)\bigl)^{[2]}=(x\dM u)^{[2]}+(x\dM v)^{[2]}+[x\dM u,x\dM v],~\forall x,u,v\in A.$

Thus, the map $(-)^{[2]}$ is well-defined. Moreover, we have
\begin{align*}
    \bigl[(x\dM u)^{[2]}, y\dM v \bigl]_{\Omega^1(A)}=&~\bigl[x^2\dM(u^{\{2\}}),y\dM v\bigl]_{\Omega^1(A)}+\bigl[x\{u,x\}\dM u,y\dM v\bigl]_{\Omega^1(A)}\\
    =&~x^2\{u^{\{2\}},y\}\dM v+x^2y\dM\{u^{\{2\}},v\}+x\{u,x\}\{u,y\}\dM v\\
    &~+yx\{\{u,x\},v\}\dM u+y\{u,x\}\{x,v\}\dM u+xy\{u,x\}\dM\{u,v\}\\
    =&~\bigl[x\dM u, x\{u,y\}\dM v + y\{x,v\}\dM u+xy\dM\{u,v\}\bigl]_{\Omega^1(A)}\\
    =&~\bigl[x\dM u,[x\dM u, y\dM v]\bigl]_{\Omega^1(A)}.
    \end{align*} Thus, by Jacobson's Theorem \ref{thm:jacobson}, $\bigl(\Omega^1(A),[-,-]_{\Omega^1(A)},(-)^{[2]}\bigl)$ is a restricted Lie algebra. It remains to show that the map $\pi^{\sharp}$ satisfies Eq. \eqref{eq:hochschild} with respect to the 2-map $(-)^{[2]}$. Let $a,x,u\in A$. We have
    \begin{align*}
        &(a\cdot x\dM u)^{[2]}+a^2\cdot(x\dM u)^{[2]}+\pi^{\sharp}(a\cdot x\dM u)(a)x\dM u\\
        =&~a^2x^2\dM(u^{[2]})+a^2x\{u,x\}\dM u+ax^2\{u,x\}\dM u+a^2\bigl(x^2\dM(u^{[2]})+x\{u,x\}\dM u \bigl)+ax\{u,a\}x\dM u\\
        =&~0.
    \end{align*} Therefore, Eq. \eqref{eq:hochschild} is satisfied on the $A$-basis of $\Omega^1(A)$. Thus, by Proposition \ref{prop:jacobson-rinehart}, $\bigl(A, \Omega^1(A),\pi^{\sharp}\bigl)$ is a restricted Lie-Rinehart algebra. Moreover, we have that
    $\pi^{\sharp}\bigl((x\dM u)^{[2]}\bigl)=\pi^{\sharp}(x\dM u)^2,$
    and since the anchor map on $\Der(A)$ is the identity of $A$, we can easily see that $\pi^{\sharp}$ is a morphism of restricted Lie algebras.
\end{proof}
\subsubsection{K\"{a}hler forms} This section aims at generalizing Eq. \eqref{eq:isomorphism deri} to multiderivations (see \cite[Section 3.2.2]{LPV}) and at extending the universal property \eqref{eq:diag-uni} to arbitrary restricted Poisson cochains. For completeness, we recall the well-known construction of K\"{a}hler forms, see \cite{LPV} for more details. Let $A$ be an \textit{unital} associative commutative $\K$-algebra. Let $k\geq 0$ and consider the $A$-modules $\Omega^k(A):=\wedge_A^k\Omega^1(A),~k\geq 1$ and $\Omega^0(A):=A.$ Moreover, we consider the graded $A$-module
$$\Omega^{\bullet}(A):=\bigoplus_{k\geq 0}\Omega^k(A),$$ whose elements are called K\"{a}hler $k$-forms. The space $\Omega^k(A)$ is generated by elements of the form $x\dM u_1\wedge\cdots\wedge \dM u_k$ as a  $\K$-vector space and by elements of the form $\dM u_1\wedge\cdots\wedge \dM u_k$ as an $A$-module, where $x,u_1,\cdots,u_k\in A.$ The map $\dM:A\rightarrow \Omega^1(A)$ extends a $\K$-linear map $\wedge^{\bullet}\dM:\wedge^{\bullet}A\rightarrow \Omega^{\bullet}(A)$ by the formula
\begin{equation}
    \wedge^{\bullet}\dM(u_1\wedge\cdots\wedge u_k):=\dM u_1\wedge\cdots\wedge\dM u_k,\quad\forall u_1,\cdots u_k\in A.
\end{equation} The K\"{a}hler $k$-forms satisfy the following universal property: for every (skew)symmetric $k$-derivation $\varphi$ of $A$, there exists an unique $A$-linear map $\hat{\varphi}:\Omega^k(A)\rightarrow A$ such that $\hat{\varphi}\circ\wedge^k\dM=\varphi,$ that is, we have
\begin{equation}\label{eq:diag-unib}
\xymatrix{
\wedge^kA \ar[rr]^(.4){\wedge^k\dM}\ar[d]_{\varphi} && \Omega^k(A)
\ar@{.>}[dll]^(.5){\exists !~\hat{\varphi}}\\
A}
\end{equation} where $A$ is seen as a module over itself using its multiplication. The map $\hat{\varphi}$ is given by
\begin{equation}\label{eq:phihat}
        \hat{\varphi}(x\dM u_1\wedge\cdots\wedge\dM u_k):=x\varphi(u_1,\cdots, u_k),\quad \forall x,u_1,\cdots, u_k\in A.
\end{equation}
As a consequence, we have a natural isomorphism of $A$-modules:
\begin{equation}\label{eq:isomorphism derivation}
\mathfrak{X}^k(A)\cong \Hom_{A}(\Omega^k(A),A)
\end{equation}
and the map $\varphi\rightarrow \hat{\varphi}$ induces a natural isomorphism 
\begin{equation}\label{eq:isomorphism derivationb}\mathfrak{X}^\bullet(A)\cong \bigoplus_{k\in \mathbb{N}}\Hom_{A}(\Omega^k(A),A).
\end{equation}

\subsubsection{The comparison theorem}
Next, we extend the diagram \eqref{eq:diag-uni} to arbitrary restricted Poisson cochains. Let $A$ be an unital Poisson algebra and let $(\varphi,\w)\in C_{\rm PA}^k(A)$ for $k\geq 2$. By \eqref{eq:diag-unib}, there exists an unique $A$-linear map $\hat{\varphi}:\Omega^k(A)\rightarrow A$ (which is given by Eq. \eqref{eq:phihat}) such that $\hat{\varphi}\circ\wedge^k\dM=\varphi$. 
\sssbegin{Proposition}\label{prop:another isomorphism}
  For $k\geq 2$, let $(\varphi,\w)\in C^k_{\rm PA}(A)$. Then, there exists an unique map $\hat{\w}:\Omega^1(A)\otimes\Omega^{k-2}(A)\rightarrow A$ such that $(\hat{\varphi},\hat{\w})\in C_{\rm LR}^k(\Omega^{k}(A);A)$ and  $\hat{\w}\circ(\dM\otimes \wedge^{k-2}\dM)=\w$, that is, the following diagram commutes:
   \begin{equation}\label{eq:diag-w}
   \xymatrix{
A\otimes \wedge^{k-2}A \ar[rr]^(.4){\quad \quad \dM\otimes \wedge^{k-2}\dM}\ar[d]_{\omega} && \Omega^1(A)\otimes \Omega^{k-2}(A)
\ar@{.>}[dll]^(.5){\exists !~\hat{\omega}}\\
A}
   \end{equation} where $A$ is seen as a module over itself using its multiplication.
\end{Proposition}
\begin{proof}
    Let $x,y,u,v,v_2,\cdots,v_{k-1}\in A.$ We define a map $\hat{\w}:\Omega^1(A)\otimes\Omega^{k-2}(A)\rightarrow A$ by the three following conditions: 
    \begin{align}\label{eq:w-hat1}
        \hat{\w}(x\dM u,y\dM v_2\wedge\dM v_3\wedge\cdots\wedge\dM v_{k-1})=&~x^2y\w(u,v_2,v_3,\cdots,v_{k-1});\\\label{eq:w-hat2}
        \hat{\w}(x\dM u+y\dM v,\dM v_2\wedge\cdots\wedge\dM v_{k-1})=&~\hat{\w}(x\dM u,\dM v_2\wedge\cdots\wedge\dM v_{k-1})\\\nonumber
        &+\hat{\w}(y\dM v,\dM v_2\wedge\cdots\wedge\dM v_{k-1})\\\nonumber
        &+\hat{\varphi}(x\dM u\wedge y\dM v\wedge\dM v_2\wedge\cdots\wedge\dM v_{k-1});\\\label{eq:w-hat3}
        \text{the map }  \dM v_2\wedge\cdots\wedge\dM v_{k-1}\mapsto  &~\hat{\w}(x\dM u,\dM v_2\wedge\cdots\wedge\dM v_{k-1}) \text{ is linear } \forall x\dM u\in \Omega^1(A).
    \end{align}
    The map $\hat{\w}$ is well-defined, since we have
    \begin{align*}
        \hat{\w}\bigl(x\dM(uv),y\dM v_2\wedge\cdots\wedge\dM v_{k-1}\bigl)&=x^2y\w\bigl(uv,v_2,\cdots, v_{k-1}\bigl)\\ 
        &=x^2y\bigl(u\w(v,v_2,\cdots, v_{k-1})+v\w(u,v_2,\cdots, v_{k-1})\\&~~+uv\varphi(u,v,v_2,\cdots, v_{k-1})\bigl)\\
        &=\hat{\w}\bigl(xu\dM v,y\dM v_2\wedge\cdots\wedge\dM v_{k-1} \bigl)+\hat{\w}\bigl(xv\dM u,y\dM v_2\wedge\cdots\wedge\dM v_{k-1} \bigl)\\
        &~~+\hat{\varphi}(xu\dM v\wedge xv\dM u\wedge y\dM v_2\wedge\cdots\wedge\dM v_{k-1}).
    \end{align*} The map $\hat{\omega}$ satisfies Eqs. \eqref{eq:res A-linear1} and \eqref{eq:res A-linear2} and uniqueness follows from the definition.
\end{proof}

\sssbegin{Theorem}\label{thm:coho-iso}
  Let $(A,\cdot,\{-,-\},(-)^{\{2\}})$ be a restricted Poisson algebra in characteristic $2$. If the module of the K\"ahler differential $\Omega^1(A)$ is free, then the cohomology complex $(\bigoplus_{n\geq0}C^{n}_{\rm PA}(A),\dM^n_{\rm PA})$ for the restricted Poisson algebra $A$ is isomorphic to the cohomology complex $\bigl(\bigoplus_{n\geq0}C^{n}_{\rm LR}(\Omega^1(A);A\bigl),\dM^n_{\rm LR})$ for the restricted Lie-Rinehart algebra $\bigl(A,\Omega^1(A), \pi^\sharp\bigl)$ described in Theorem \ref{thm:RLRfromPoisson}.
\end{Theorem}
\begin{proof}
By \eqref{eq:isomorphism derivation} and Proposition \ref{prop:another isomorphism}, the map 
\begin{equation}
   \Phi:\, \bigoplus_{n\geq0}C^{n}_{\rm PA}(A) \rightarrow\bigoplus_{n\geq0}C^{n}_{\rm LR}(\Omega^1(A);A\bigl),\quad (\varphi,\omega)\in C^{n}_{\rm PA}(A)\mapsto (\hat{\varphi},\hat{\omega})\in C^{n}_{\rm LR}(\Omega^1(A);A\bigl).
\end{equation} defines an isomorphism between the graded spaces $\bigoplus_{n\geq0}C^{n}_{\rm PA}(A)\cong \bigoplus_{n\geq0}C^{n}_{\rm LR}(\Omega^1(A);A\bigl)$.

Let $(\varphi,\w)\in C^{n}_{\rm PA}(A)$.  As in the ordinary case, we have
$\Phi (\dM^n_{\rm PA}(\varphi))=\dM^n_{\rm LR}(\hat{\varphi})$, see \cite[Proposition 4.5]{LPV}.

Next, we show that $\Phi (\delta^n(\omega))=\delta^n(\hat{\omega})$ holds. Let $u,v_2,\cdots,v_n\in A$. We have
\begin{eqnarray*}
& & \Phi \bigl(\delta^n(\omega)\bigl)(\dM u,\dM v_2\wedge\cdots\wedge \dM v_n)  \\
&=& \delta^n(\omega)(u,v_2,\cdots,v_n)   \\
&=&\{u,\varphi(u,v_2,\cdots,v_n)\}+\sum_{i=2}^{n}\{v_i,\omega(u,v_2,\cdots,\hat{v_i},\cdots,v_n)\}+\varphi(u^{\{2\}},v_2,\cdots,v_n)\\	
&&+\sum_{i=2}^{n}\varphi\left(\{u,v_i\},u,v_2,\cdots,\hat{v_i},\cdots,v_n \right)+\sum_{2\leq i<j\leq n}\omega\left(u,\{v_i,v_j\},v_2,\cdots,\hat{v_i},\cdots,\hat{v_j},\cdots,v_n  \right)\\
    &=&\pi^\sharp(du)\hat{\varphi}(\dM u\wedge\dM v_2\wedge\cdots\wedge\dM v_n)+\sum_{i=2}^{n}\pi^\sharp(\dM v_i)\hat{\omega}(\dM u,\dM v_2\wedge\cdots\wedge\widehat{\dM v_i}\wedge\cdots\wedge\dM v_n)\\
    &&+\hat{\varphi}\bigl((\dM u)^{[2]}\wedge\dM v_2\wedge\cdots\wedge\dM v_n\bigl)+\sum_{i=2}^{n}\hat{\varphi}\bigl([\dM u,\dM v_i]_{\Omega^1(A)}\wedge\dM u\wedge\dM v_2\wedge\cdots\wedge\widehat{\dM v_i}\wedge\cdots\wedge\dM v_n \bigl)\\
    &&+\sum_{2\leq i<j\leq n}\hat{\omega}\left(\dM u,[\dM v_i,\dM v_j]_{\Omega^1(A)}\wedge\dM v_2\wedge\cdots\wedge\widehat{dv_i}\wedge\cdots\wedge\widehat{\dM v_j}\wedge\cdots\wedge\dM v_n  \right)\\
    &=&\delta^n(\hat{\omega})(\dM u,\dM v_2\wedge\cdots\wedge\dM v_n).
\end{eqnarray*}
Moreover, since $\Phi (\dM^n_{\rm PA}(\omega))$ and $\dM^n_{\rm res}(\hat{\omega})$ are restricted cochains satisfying \eqref{eq:res A-linear1}-\eqref{eq:res A-linear2}, and since $\hat{\varphi}$ and $\hat{\w}$ satisfy Eqs. \eqref{eq:phihat},\eqref{eq:w-hat1},\eqref{eq:w-hat2} and \eqref{eq:w-hat3}, this relation extends to any arbitrary elements $\alpha,\beta_2,\cdots,\beta_n\in \Omega^1(A)$:
$$\Phi \bigl(\delta^n(\omega)\bigl)(\alpha,\beta_2,\cdots,\beta_n)=\delta^n(\hat{\omega})(\alpha,\beta_2,\cdots,\beta_n).$$
Therefore, the cohomology complex $(\bigoplus_{n\geq0}C^{n}_{\rm PA}(A),\dM^n_{\rm PA})$ for the restricted Poisson algebra $A$ is isomorphic to the cohomology complex $\bigl(\bigoplus_{n\geq0}C^{n}_{\rm LR}(\Omega^1(A);A\bigl),\dM^n_{\rm LR})$ for the restricted Lie-Rinehart algebra $\bigl(A,\Omega^1(A), \pi^\sharp\bigl)$.
\end{proof}

\subsection{Formal deformations of restricted Poisson algebras}\label{sec:defopoisson}

Let $t$ be a formal parameter, and 
denote by $\K^t$ the ring $\K[[t]]$ of formal power series in $t$ with coefficients in $\K$. Let $(A,\cdot)$ be a commutative  associative algebra. We denote by $A^t$ the $\K^t$-algebra of formal power series with coefficients in $A$. The product in $A^t$ is given by
\begin{equation}
\bigl(\sum_{i\geq0}{x_it^i}\bigl)\cdot \bigl(\sum_{j\geq0}y_jt^j\bigl)=\sum_{i,j\geq0}x_iy_jt^{i+j},\quad\forall x_i,y_i\in A.
\end{equation}
This product is associative, $\K^t$-bilinear and reduces to the original product $\cdot$ of $A$ when $t=0$.

Let $\bigl(A,\cdot,\{-,-\},(-)^{\{2\}}\bigl)$ be a restricted Poisson algebra. Consider two maps
\begin{align}
    \mu_\star: A^t\times A^t\rightarrow A^t,\quad \text{ and }\quad \w_\star:A^t\rightarrow A^t.
\end{align}

A tuple $\bigl(A^t,\cdot,\mu_\star,\w_\star\bigl)$ is called a \textit{formal deformation} of the restricted Poisson algebra $\bigl(A,\cdot,\{-,-\},(-)^{\{2\}}\bigl)$ if $(A^t,\mu_\star,\w_\star)$ is a restricted Lie algebra and if in addition, we have
\begin{align}\label{eq:defo-poisson1}
    \mu_\star(x\cdot y,z)&=x\cdot\mu_\star(y,z)+y\cdot\mu_\star(x,z), ~\quad\forall x,y,z\in A;\\\label{eq:defo-poisson2}
    \w_\star(x\cdot y)&=x^2\cdot\w_\star(y)+y^2\cdot\w_\star(x)+x\cdot y\cdot\mu_\star(x,y),\quad\forall x,y\in A,
\end{align} that is, $\bigl(A^t,\cdot,\mu_\star,\w_\star\bigl)$ is a restricted Poisson algebra.

Two formal deformations $\bigl(A^t,\cdot,\mu_\star,\w_\star\bigl)$ and $\bigl(A^t,\cdot,\mu'_\star,\w'_\star\bigl)$ of $\bigl(A,\cdot,\{-,-\},(-)^{\{2\}}\bigl)$ are called \textit{equivalent} if there exists an isomorphism of restricted Poisson algebras \begin{equation}
\Psi:\bigl(A^t,\cdot,\mu_\star,\w_\star\bigl)\rightarrow\bigl(A^t,\cdot,\mu'_\star,\w'_\star\bigl),
\end{equation}
of the form
$$\Psi=\id+\sum_{i\geq1}t^i\psi_i,$$ where $\psi_i:A\rightarrow A $ are linear maps.

Formal deformations are usually constructed term by term, which is done by considering $k$-th order deformations. For $k\geq 0$, we consider the ring $\K^t_k=\K^t/<t^{k+1}>$. Starting from an associative commutative associative algebra $A$, we consider similarly $A^t_k=A^t/<t^{k+1}>$, which can be seen as an $\K_k^t$-algebra.

Consider two maps given by
\begin{align}
    \mu_{(k)}:=\{-,-\}+\sum_{i\geq 1}^kt^i\mu_i,\quad \text{ and }\quad \w_{(k)}:=(-)^{\{2\}}+\sum_{i\geq 1}^kt^i\w_i,
\end{align} where $(\mu_i,\w_i)\in C^2_{\rm PA}(A)~\forall i\geq 1.$ A tuple $\bigl(A^t_k,\cdot,\mu_{(k)},\w_{(k)}\bigl)$ is called a \textit{formal deformation of order} $k$ of the restricted Poisson algebra $\bigl(A,\cdot,\{-,-\},(-)^{\{2\}}\bigl)$ if $(A^t_k,\mu_{(k)},\w_{(k)})$ is a restricted Lie algebra. Since $(\mu_i,\w_i)\in C^2_{\rm PA}(A)~\forall i\geq 1,$ the maps $\mu_{(k)}$ and $\w_{(k)}$ satisfy Eqs. \eqref{eq:defo-poisson1} and \eqref{eq:defo-poisson2}. Therefore, $\bigl(A_k^t,\cdot,\mu_{(k)},\w_{(k)}\bigl)$ is a restricted Poisson algebra.

The formal deformations of order $k=1$ are called \textit{infinitesimal}. 

\sssbegin{Lemma}
Let $k\geq 1$ and let $\bigl(A_k^t,\cdot,\mu_{(k)},\w_{(k)}\bigl)$ be a formal deformation of order $k$ of a restricted Poisson algebra $\bigl(A,\cdot,\{-,-\},(-)^{\{2\}}\bigl)$. Then, we have $(\mu_1,\w_1)\in Z^2_{\rm PA}(A)$.
\end{Lemma}
\begin{proof}
    It is clear that $(\mu_1,\w_1)\in C^2_{\rm PA}(A)$. It remains to show the restricted 2-cocycle condition, which is covered by \cite[Proposition 4.1.3]{EM}.
\end{proof}
\sssbegin{Proposition}
    Let $k\geq 1$, and let $\bigl(A^t_k,\cdot,\mu_{k},\w_{k}\bigl)$ and $\bigl(A^t_k,\cdot,\mu'_{k},\w'_{k}\bigl)$ be two equivalent deformations of $\bigl(A,\cdot,\{-,-\},(-)^{\{2\}}\bigl)$. Then,
    $(\mu_1,\w_1)\equiv (\mu_1',\w_1')\in \mathrm{H}^2_{\rm PA}(A).$
\end{Proposition}
\begin{proof}
    Suppose that the equivalence is realized by $\Psi=\id+\sum_{i\geq1}t^i\psi_i,$ where $\psi_i:A\rightarrow A $ are linear maps. For all $x,y\in A$, we have
    \begin{align*}
        \Psi(xy)=\Psi(x)\Psi(y)&=\bigl(x+t\psi_1(x)\bigl)\bigl(y+t\psi_1(y)\bigl)\quad \mod t^2\\
        &=xy+t\bigl(\psi_1(x)y+x\psi_1(y)\bigl)\quad \mod t^2.
    \end{align*}
    Since $\Psi(xy)=xy+t\psi_1(xy)\,\mod t^2$, we obtain that $\psi_1(xy)=\psi_1(x)y+x\psi_1(y)$. It follows that  $\psi_1\in C^1_{\rm PA}(A)$. Then, by \cite[Lemma 4.15]{EM}, we have $(\mu_1,\w_1)=(\mu_1',\w_1')+\dM^1_{\rm PA}\psi_1.$
\end{proof}
\sssbegin{Corollary}\label{cor:equiv-defo}
   Let $\bigl(A,\cdot,\{-,-\},(-)^{\{2\}}\bigl)$ be a restricted Poisson algebra. Then, the second cohomology space $\mathrm{H}^2_{\rm PA}(A)$ classifies up to equivalence the infinitesimal deformations of $\bigl(A,\cdot,\{-,-\},(-)^{\{2\}}\bigl)$.
\end{Corollary}

In what follows, we discuss the problem of extending a deformation of order $k$ to order $k+1$, for $k\geq 1$.   Let $\bigl(A^t_k,\cdot,\mu_{(k)},\w_{(k)}\bigl)$ be a formal deformation of order $k$ of a restricted Poisson algebra $\bigl(A,\cdot,\{-,-\},(-)^{\{2\}}\bigl)$. Let $x,y,z\in A$, and let us define the maps
\begin{align}
    \obs^{(1)}_{k+1}(x,y,z)&:=\sum_{i=1}^k\Bigl(\mu_i(x,\mu_{k+1-i}(y,z))+\mu_i(y,\mu_{k+1-i}(z,x))+\mu_i(z,\mu_{k+1-i}(x,y))   \Bigl); \\
    \obs^{(2)}_{k+1}(x,y)&:=\sum_{i=1}^k\Bigl(\mu_i\bigl(y,\w_{k+1-i}(x)\bigl)+\mu_i\bigl(x,\mu_{k+1-i}(x,y)\bigl)    \Bigl).
\end{align}

\sssbegin{Lemma}
    The pair $\bigl(\obs^{(1)}_{k+1},\obs^{(2)}_{k+1}\bigl)$ is a 3-cochain of the restricted Poisson cohomology, that is, $\bigl(\obs^{(1)}_{k+1},\obs^{(2)}_{k+1}\bigl)\in C^3_{\rm PA}(A)$.
\end{Lemma}

\begin{proof}
    By \cite[Lemma 4.20]{EM}, we have $\bigl(\obs^{(1)}_{k+1},\obs^{(2)}_{k+1}\bigl)\in C^3_{\rm res}(A;A)$. By a direct computation, we have
    \begin{equation}
        \obs^{(1)}_{k+1}(uv,y,z)=u\cdot\obs^{(1)}_{k+1}(v,y,z)+v\cdot\obs^{(1)}_{k+1}(u,y,z),~\forall u,v,y,z\in A.
    \end{equation} Moreover, for all $x,y,z\in A$, we have
\footnotesize{\begin{align*}
    &\obs^{(2)}_{k+1}(xy,z)\\=&~\sum_{i=1}^k\Bigl(\mu_i\bigl(z,\w_{k+1-i}(xy)\bigl)+\mu_i\bigl(xy,\mu_{k+1-i}(xy,z)\bigl) \\
    =&~\sum_{i=1}^k\Bigl(x^2\mu_i(z,\w_{k+1-i}(y))+y^2\mu_i(z,\w_{k+1-i}(x)) +xy\cdot\mu_i(z,\mu_{k+1-i}(x,y))+\mu_{k+1-i}(x,y)\mu_i(z,xy) \Bigl)\\
    &~+\sum_{i=1}^k\Bigl( x^2\mu_i(y,\mu_{k+1-i}(y,z))+x\mu_{k+1-i}(y,z)\mu_i(y,x)+xy\cdot\mu_i(y,\mu_{k+1-i}(x,z))     \Bigl)\\
    &~+\sum_{i=1}^k\Bigl( y^2\mu_i(x,\mu_{k+1-i}(x,z))+y\mu_{k+1-i}(x,z)\mu_i(y,x)+xy\cdot\mu_i(x,\mu_{k+1-i}(y,z))     \Bigl)\\
    =&~x^2\obs^{(2)}_{k+1}(y,z)+y^2\obs^{(2)}_{k+1}(x,z)+xy\cdot\obs^{(1)}_{k+1}(x,y,z)\\
    &~+\sum_{i=1}^k\Bigl(x\mu_{k+1-i}(x,y)\mu_i(z,y)+y\mu_{k+1-i}(x,y)\mu_i(z,x)+x\mu_{k+1-i}(y,z)\mu_i(y,x)+y\mu_{k+1-i}(x,z)\mu_i(y,x)  \Bigl)\\
    =&~x^2\obs^{(2)}_{k+1}(y,z)+y^2\obs^{(2)}_{k+1}(x,z)+xy\cdot\obs^{(1)}_{k+1}(x,y,z),
\end{align*}}\normalsize{} since the remaining sum vanishes. It follows that $\bigl(\obs^{(1)}_{k+1},\obs^{(2)}_{k+1}\bigl)\in C^3_{\rm PA}(A)$.
\end{proof}

\sssbegin{Proposition}\label{prop:obs}
     Let $\bigl(A^t_k,\cdot,\mu_{(k)},\w_{(k)}\bigl)$ be a formal deformation of order $k$ of a restricted Poisson algebra $\bigl(A,\cdot,\{-,-\},(-)^{\{2\}}\bigl)$. Let $(\mu_{k+1},\w_{k+1})\in C^2_{\rm PA}$. Then, $\bigl(A^t_{k+1},\cdot,\mu_{(k+1)},\w_{(k+1)}\bigl)$ is a formal deformation of order $k+1$ of a restricted Poisson algebra $\bigl(A,\cdot,\{-,-\},(-)^{\{2\}}\bigl)$ if and only if
     $	\bigl( \obs^{(1)}_{k+1}, \obs^{(2)}_{k+1}\bigl)=\dM_{\rm PA}^2\left(\mu_{k+1},~\w_{k+1} \right).        $
     \end{Proposition}

     \begin{proof}
Suppose that $\bigl(A^t_{k+1},\cdot,\mu_{(k+1)},\w_{(k+1)}\bigl)$ is a formal deformation of order $k+1$.  By \cite[Proposition 4.21]{EM}, we have       $\bigl( \obs^{(1)}_{k+1}, \obs^{(2)}_{k+1}\bigl)=\dM_{\rm PA}^2\left(\mu_{k+1},~\w_{k+1} \right).$ Since we are working over $A^t_{k+1}$, this is actually an equivalence. 
     \end{proof}

\subsection{Examples of restricted Poisson algebras}\label{sec:ex}

In this section, we present some examples of restricted Poisson algebras and compute their cohomology. The survey paper \cite{KKP} contains a classification of Poisson algebras of small dimension (in characteristic 0). The examples we borrow from \cite{KKP} are also valid in characteristic 2. 
\subsubsection{Heisenberg algebra of dimension 3} Consider the 3-dimensional restricted Heisenberg Lie algebra $\fh$ spanned by elements $e_1,e_2,e_3$ with the bracket $[e_1,e_2]=e_3$, and the 2-map $(\cdot)^{[2]}= 0$. We can endow this Lie algebra with an associative commutative product given by $e_1e_2=e_3$.  Therefore, $\fh$ is a restricted Poisson algebra.

We shall compute the first and second cohomology groups $\mathrm{H}_{\rm PA}^1(\fh)$ and $\mathrm{H}_{\rm PA}^2(\fh)$. A basis of the space $\mathfrak{X}^1(\fh)$ is given by
\begin{equation*}
 \begin{array}{lllllllll}
        \psi_1&=&e_1\otimes e_1^*+e_2\otimes e_2^*,&\psi_2&=&e_1\otimes e_1^*+e_3\otimes e_3^*, &\psi_3&=&e_1\otimes e_2^*,\\
        \psi_4&=&e_2\otimes e_1^*,&\psi_5&=&e_3\otimes e_1^*, &\psi_6&=&e_3\otimes e_1^*.
\end{array}
\end{equation*}
A direct computation shows that the maps $\psi_1,\psi_2,\psi_5,\psi_6$ are restricted Poisson 1-cocycles while $\psi_5$ and $\psi_6$ are restricted Poisson 1-coboundaries. Therefore,
\begin{equation}
    \mathrm{H}_{\rm PA}^1(\fh)=\Span\{\psi_1,\psi_2\}.
\end{equation}

Similarly, we have 
\begin{equation}
    \mathfrak{X}^1(\fh)\cap Z^2_{\ce}(\fh;\fh)=\bigl\{\varphi_1=e_1\otimes (e_1^*\wedge e_2^*)+e_3\otimes (e_2^*\wedge e_3^*); \varphi_2=e_2\otimes (e_1^*\wedge e_2^*)+e_3\otimes (e_1^*\wedge e_3^*) \bigl\}. 
\end{equation} Moreover, a map $\w_1:\fh\rightarrow\fh$ satisfies Eq. \eqref{eq:res PA1} with respect to $\varphi_1$ (on the basis) if and only if $\w_1(e_3)=0.$ In addition,  the restricted cocycle condition $\delta^2\w_1=0$ leads to $\w_1(e_1)=\lambda e_3,~\w_1(e_2)=\mu e_3,$ for some $\lambda,\mu\in \K$. But, $\w_1$ also satisfies Eq. \eqref{eq:res PA1} with respect to $\varphi\equiv 0$, and we have $(0,\w_1)=\dM^1_{\rm PA}(\lambda e_2\otimes e_1^*+\mu e_1\otimes e_2^*)$. Therefore, we have 
\[
(\varphi_1, \omega_1)=(0,\omega_1)+(\varphi_1,0) \cong (\varphi_1,0).
\]
Similarly, we obtain another nontrivial restricted Poisson 2-cocycle $(\varphi_2,0)$. Therefore, we have
\begin{equation}
    \mathrm{H}_{\rm PA}^2(\fh)=\bigl\{ (\varphi_1,0),(\varphi_2,0) \bigl\}.
\end{equation}
It is worth noticing that $\mathrm{H}^1_{\mathrm{CE}}(\mathfrak{h};\mathfrak{h})$ is 4-dimensional while $\mathrm{H}^2_{\mathrm{CE}}(\mathfrak{h};\fh)$ is 2-dimensional, see \cite{BGLS}
\subsubsection{A rigid restricted Poisson algebra of dimension 3} \label{ex:rigiddim3}  Consider the 3-dimensional restricted Lie algebra $L$ spanned by elements $e_1,e_2,e_3$ with the bracket $[e_2,e_3]=e_2$, and the the 2-map $e_3^{[2]}= e_3$. We can endow this Lie algebra with an associative commutative product given by $e_1e_1=e_1$.  Therefore, $L$ is a restricted Poisson algebra.
We shall compute the first and second cohomology groups $\mathrm{H}_{\rm PA}^1(L)$ and $\mathrm{H}_{\rm PA}^2(L)$. A basis of the space $Z_{\rm PA}^1(L)$ is given by $\bigl\{e_2\otimes e_2^*,e_2\otimes e_3^*\bigl\}.$ But since $e_2\otimes e_2^*=\ad_{e_3} \text{ and } e_2\otimes e_3^*=\ad_{e_2}. $
It follows that 
\begin{equation}
    \HH^1_{\rm PA}(L)=0.
\end{equation}
Similarly, we have 
\begin{equation}
    \mathfrak{X}^1(L)\cap Z^2_{\ce}(L;L)=\bigl\{\varphi_1=e_2\otimes (e_2^*\wedge e_3^*), \varphi_2=e_3\otimes (e_2^*\wedge e_3^*) \bigl\}. 
\end{equation}
A careful computation shows that
\begin{equation}
    Z^2_{\rm PA}(L)=\bigl\{(\varphi_1,\w_1), (\varphi_2,\w_2) \bigl\},
\end{equation} where $\w_1(e_3)=e_3$ and $\w_2(e_2)=e_2$. But since  
\[
(\varphi_1,\w_1)=\dM_{\rm PA}^1(e_3\otimes e_3^*) \text{ and } (\varphi_2,\w_2)=\dM_{\rm PA}^1(e_3\otimes e_2^*),
\]
it follows that 
\begin{equation}
    \HH^2_{\rm PA}(L)=0.
\end{equation}
It is worth mentioning that $\dim\HH^1_{\ce}(L;L)=2 \text{ while } \dim\HH^2_{\ce}(L;L)=1.$

\subsubsection{A non-rigid restricted Poisson algebra of dimension 3} We consider the same restricted Lie algebra $L$ of dimension 3 as in Example \ref{ex:rigiddim3}, but endowed with the associative commutative multiplication $e_1e_1=e_1,~e_1e_2=e_2,$ and $e_1e_3=e_3$. In that case, we have
\begin{equation}
     \HH^1_{\rm PA}(L)=0, \text{ and } \HH^2_{\rm PA}(L)=\Span\{(0,\w)\}, \text{ where } \w(e_2)=e_1.
\end{equation}

\subsubsection{Analogs of the classical Poisson algebras $\mathfrak{po}_{\Pi}(2n,\underline{N})$} Let $k\geq 1$. The associative commutative algebra of divided powers in $k$ variables $x:=(x_1,\cdots,x_k)$ is defined for $\underline{N}:=(n_1,\cdots,n_k)$, where $n_s\geq 0,~\forall~  1\leq s \leq k$, by (see, e.g,  \cite{Ko, SF})
\begin{align}
    \K(x;\underline{N}):=\Span\bigl\{x^{(\underline{i})}:=x_1^{(i_1)}\cdots x_k^{(i_k)},~(\underline{i})=(i_1,\cdots,i_k),~0\leq i_s\leq p^{n_s}-1  \bigl\}.
\end{align}
The multiplication is given by
\begin{equation}\label{eq:divpwr}
x^{(\underline{i})}x^{(\underline{j})}=\binom{\underline{i}+\underline{j}}{\underline{i}}x^{(\underline{i}+\underline{j})},~\text{where}~\binom{\underline{i}+\underline{j}}{\underline{i}}:=\prod_{s=1}^k\binom{i_s+j_s}{i_s}.
\end{equation}

For $n\geq 1$, we define (see, \cite{LeD})
 \begin{equation}
    \po_{\Pi}(2n,\underline{N}):=\{f\in\K[p_i,q_i;\underline{N}],~1\leq i\leq n\},
 \end{equation}
endowed with the Lie bracket 
\begin{equation}
    \{f,g\}_{\Pi}:=\sum_{i=1}^n\left(\frac{\del f}{\del p_i}\frac{\del g}{\del q_i}-\frac{\del f}{\del q_i}\frac{\del g}{\del p_i}\right ).
  \end{equation}
With the divided powers multiplication (see Eq. \eqref{eq:divpwr}
), $\po(2n,\underline{N})$ is a Poisson algebra. The Lie algebra $\po(2n,\underline{N})$ is restricted if and only if $\underline{N}=\underline{1}:=(1,\cdots,1)$. In that case, the 2-map is given by 
\begin{align*}
1^{\{2\}}=\mu 1,~p_i^{\{2\}}=\lambda_{i,0}1,~q_i^{\{2\}}&=\lambda_{0,i}1,~(p_iq_j)^{\{2\}}=\delta_{i,j}p_iq_j+\lambda_{i,j} 1,\\ (p_{i_1}\cdots p_{i_s}q_{j_1}\cdots q_{j_t})^{\{2\}}&=\lambda_{i_1,\cdots, i_s, j_1,\cdots, j_t}1,\quad\forall s+t\geq 3,
\end{align*}
where the $\lambda$'s and $\mu$ belong to $\K$ and where $\delta$ is the Kronecker symbol. Because of the presence of the center, Skryabin \cite{Sk}, introduced the notion of {\it normalized} $2$-map on $\po(2n,\underline{1})$ by assuming that $1^{\{2\}}=0$ and $f^{\{2\}} \in \mathfrak{m}^2$ for all $f \in \mathfrak{m}^2$, where $\mathfrak m$ is the maximal ideal of $\po(2n,\underline{1})$ as an associative algebra. The Lemma \ref{Skry} below shows that Condition \eqref{eq:resPA-2map} forces the $2$-map introduced here to be  normalized. 

\sssbegin{Lemma}\label{Skry} For all $n\geq 1$, the restricted Lie algebra $\po_{\Pi}(2n,\underline{1})$ is a restricted Poisson algebra if and only if $\mu=\lambda_{i_1,\cdots, i_s, j_1,\cdots, j_t}=0,\quad\forall s,t\geq 0.$
\end{Lemma}
\begin{proof} Straightforward.
\end{proof}
We have
\begin{align}
    \HH^1_{\ce}\bigl(\po_{\Pi}(2,\underline{1}),\po_{\Pi}(2,\underline{1})\bigl)=\Span\bigl\{ 1\otimes \left(p_1\, q_1\right)^*, p_1\otimes q_1^*,~q_1\otimes p_1^*,~ p_1\otimes p_1^*+1\otimes 1^*\bigl\}
\end{align}
We will provide an example here in the case of $\po_{\Pi}(2,\underline{1})\cong \mathfrak{gl}(2).$
\sssbegin{Proposition} 
    We have $\HH^1_{\PA}\bigl(\po_{\Pi}(2,\underline{1})\bigl)=0$ and $\HH^2_{\PA}\bigl(\po_{\Pi}(2,\underline{1})\bigl)=\Span\bigl\{(0,\w)\bigl\},$ where $$\w(1)=1,~\w(p_1)=\w(q_1)=\w(p_1q_1)=0.$$
\end{Proposition}

\begin{proof}
  Computing the dimensions of $Z^1_{\PA}\bigl(\po_{\Pi}(2,\underline{1})\bigl)$ and $B^1_{\PA}\bigl(\po_{\Pi}(2,\underline{1})\bigl)$ using SuperLie leads to $\HH^1_{\PA}\bigl(\po_{\Pi}(2,\underline{1})\bigl)=0$. Next, we have that $\mathfrak{X}^2\bigl(\po_{\Pi}(2,\underline{1})\bigl)\cap Z^2_{\ce}\bigl(\po_{\Pi}(2,\underline{1});\po_{\Pi}(2,\underline{1})\bigl)$ is one dimensional and spanned by $\varphi=(p_1q_1)\otimes(p_1^*\wedge q_1^*).$ Direct computations show that there is no map $\w:\po_{\Pi}(2,\underline{1})\rightarrow \po_{\Pi}(2,\underline{1})$ satisfying \eqref{eq:res PA1} such that $\delta^2\w=0$. Therefore, the map $\varphi$ does not produce any restricted Poisson cocycle.
  
  It remains to investigate cocycles of the form $(0,\w)$. A pair $(0,\w)$ is a restricted cocycle satisfying \eqref{eq:res PA1} only if there exists scalars $\lambda,\mu,\gamma$ such that,
  $$\w(p_1q_1)=0;~\w(p_1)=\lambda 1~;\w(q_1)=\mu 1~;\w(1)=\gamma 1.$$ But we have $\delta^1(\mu p_1\otimes q_1+\lambda q_1\otimes p_1)=(\lambda+\mu)1,$ thus the only non-trivial 2-cocycles are given by pairs $(0,\w)$, with $\w(1)=\gamma1,~\gamma\in \K$.
\end{proof}
It is worth noticing that $\dim\HH^2_{\ce}\bigl(\po_{\Pi}(2,\underline{1})\bigl)=6.$
\sssbegin{Proposition} 
    We have $\HH^1_{\PA}\bigl(\po_{\Pi}(4,\underline{1})\bigl)=0$.
\end{Proposition}
\begin{proof}
Long computation assisted by SuperLie (\cite{G}).
\end{proof}

\noindent\textbf{Remark.} In characteristic 2, there is another version of Poisson algebra, see \cite{LeD, BGLLS}.
For $n\geq 1$, we let
 \begin{equation}
    \po_{\rm I}(n,\underline{N}):=\{f\in\K[z_i;\underline{N}],~1\leq i\leq n\},
 \end{equation}
with the bracket 
\begin{equation}
    \{f,g\}_{\rm I}:=\sum_{i=1}^{n}\frac{\del f}{\del z_i}\frac{\del g}{\del z_i}.
  \end{equation}
This is not a Lie bracket. Indeed, we have $\{z_i,z_i\}=1,~\forall 1\leq i \leq n.$

We then consider $\fh_{\rm I}(n,\underline{N}):=\po_{\rm I}(n,\underline{N})/\Span\{1\}$ with the same bracket. Although $\bigl(\fh_{\rm I}(2n,\underline{N}),\{\cdot,\cdot\}_{\rm I}\bigl)$ is now a Lie algebra, nevertheless it is not restricted.\\

\noindent\textbf{Acknowledgment.} We thank V. Petrogradsky for pointing out the reference \cite{PS}.


\end{document}